\documentclass[preprint,12pt]{elsarticle}
\journal{ }

\usepackage[utf8]{inputenc}
\usepackage{amsmath, amsxtra, amsfonts,amscd,amssymb,amsthm,mathtools}
\usepackage{algorithm,algorithmic}
\usepackage{color}
\usepackage{graphicx,subfigure}
\usepackage{epstopdf}
\usepackage{url}
\usepackage[left=1in,top=1in,right=1in,bottom=1in,letterpaper]{geometry}
\usepackage{setspace}
\usepackage{multirow,makecell,booktabs}
\usepackage[table,xcdraw]{xcolor}

\newcommand{\vx}{\mathbf{x}}
\newcommand{\vy}{\mathbf{y}}
\newcommand{\mfd}{\mathcal{M}}
\newcommand{\mtc}{g}
\newcommand{\ctl}{\mathbf{v}}
\newcommand{\tgsp}{\mathrm{T}\mfd}
\newcommand{\tgspx}{\mathrm{T}_{\vx}\mfd}
\newcommand{\vm}{\mathbf{m}}
\newcommand{\cstr}{\mathcal{C}}
\newcommand{\prob}{\mathcal{P}}

\newcommand{\vfunc}{\phi}
\newcommand{\Jfunc}{J}
\newcommand{\costfuncctl}{L}
\newcommand{\hmtctl}{H}
\newcommand{\costfuncevo}{F}
\newcommand{\costfuncterm}{F_{T}}
\newcommand{\vp}{\mathbf{p}}
\newcommand{\vq}{\mathbf{q}}
\newcommand{\vu}{\mathbf{u}}

\newcommand{\kernel}{K}

\newcommand{\Costevo}{\mathcal{F}}
\newcommand{\Costterm}{\mathcal{F}_{T}}
\newcommand{\Cost}{\mathcal{Y}}

\newcommand{\tk}{{t_k}}
\newcommand{\tkp}{{t_{k+\half}}}
\newcommand{\tkm}{{t_{k-\half}}}
\newcommand{\tkmm}{t_{k-1}}
\newcommand{\tn}{{t_n}}

\newcommand{\disctmfd}{\widetilde{\mfd}}
\newcommand{\discttgsp}{\mathrm{T}\disctmfd}
\newcommand{\vtc}{V}
\newcommand{\trg}{T}
\newcommand{\vtci}{{\vtc_i}}
\newcommand{\trgj}{{\trg_j}}
\newcommand{\numvtc}{h}
\newcommand{\numtrg}{s}
\newcommand{\area}{A}
\newcommand{\areavtc}{\area_{\vtc}}
\newcommand{\areavtci}{\area_{\vtci}}
\newcommand{\areatrg}{\area_{\trg}}
\newcommand{\areatrgj}{\area_{\trgj}}

\newcommand{\sumtrg}{\sum_{j:\vtci\in\trgj}}
\newcommand{\disctmtc}{G}
\newcommand{\Psibar}{\overline{\Psi}}

\newcommand{\grad}{\nabla_{\disctmfd}}
\newcommand{\divg}{\nabla_{\disctmfd}\cdot}

\newcommand{\Zero}{\mathbf{0}}
\newcommand{\Rho}{P}
\newcommand{\Rhobar}{\overline{\Rho}}

\newcommand{\Dt}{\widetilde{\partial}_t}

\newcommand{\disctcstr}{\widetilde{\cstr}}
\newcommand{\disctcostfuncctl}{\widetilde{\costfuncctl}}

\newcommand{\disctCostevo}{\widetilde{\Costevo}}
\newcommand{\disctCostterm}{\widetilde{\Costterm}}
\newcommand{\disctCost}{\widetilde{\Cost}}
\newcommand{\discthmt}{\widetilde{\hmtctl}}
\newcommand{\disctkernel}{\widetilde{\kernel}}
\newcommand{\vtctotrg}{w}
\newcommand{\Vtctotrg}{W}

\newcommand{\nit}{{(l)}}
\newcommand{\nitp}{{(l+\half)}}
\newcommand{\nitpp}{{(l+1)}}

\newcommand{\half}{\frac{1}{2}}
\newcommand{\third}{\frac{1}{3}}
\newcommand{\bbR}{\mathbb{R}}
\newcommand{\Lag}{\mathcal{A}}

\newcommand{\dd}{\mathrm{d}}
\newcommand{\ddx}{\dd_\mfd \vx}
\newcommand{\ddy}{\dd_\mfd \vy}

\DeclareMathOperator{\trace}{tr}
\DeclareMathOperator{\diag}{diag}

\DeclareMathOperator{\proj}{proj}
\DeclareMathOperator*{\argmin}{argmin}

\DeclareMathOperator*{\spn}{span}

\newtheorem{theorem}{Theorem}[section]

\theoremstyle{definition}
\newtheorem{definition}[theorem]{Definition}
\newtheorem{example}[theorem]{Example}

\theoremstyle{remark}
\newtheorem{remark}[theorem]{Remark}

\begin{document}

\begin{frontmatter}

\title{Computational Mean-field Games on Manifolds}
\author[1]{Jiajia Yu\fnref{fn1}}
\ead{yuj12@rpi.edu}

\author[1]{Rongjie Lai\corref{cor}\fnref{fn1}}
\ead{lair@rpi.edu}

\author[2]{Wuchen Li\fnref{fn2}}
\ead{wuchen@mailbox.sc.edu}

\author[3]{Stanley Osher\fnref{fn3}}
\ead{sjo@math.ucla.edu}

\cortext[cor]{Corresponding author}
\fntext[fn1]{J. Yu and R. Lai's work are supported in part by an NSF Career Award DMS–1752934 and NSF DMS-2134168.}
\fntext[fn2]{W. Li's work is supported in part by AFOSR MURI FP 9550-18-1-502 and NSF RTG: 2038080.}
\fntext[fn3]{S. Osher's work is supported in part by AFOSR MURI FP 9550-18-1-502, and ONR
grants:  N00014-20-1-2093, and N00014-20-1-2787.}

\affiliation[1]{organization={Department of Mathematics, Rensselaer Polytechnic Institute},
                city={Troy},
                postcode={NY 12180},
                country={USA}}
\affiliation[2]{organization={Department of Mathematics, University of South Carolina},
                city={Columbia},
                postcode={SC 29208},
                country={USA}}
\affiliation[3]{organization={Department of Mathematics, University of California, Los Angeles},
                city={Los Angeles},
                postcode={CA 90095},
                country={USA}}                

\begin{abstract}
Conventional Mean-field games/control study the behavior of a large number of rational agents moving in the Euclidean spaces. In this work, we explore the mean-field games on Riemannian manifolds. We formulate the mean-field game Nash Equilibrium on manifolds. We also establish the equivalence between the PDE system and the optimality conditions of the associated variational form on manifolds. Based on triangular mesh representation of two-dimensional manifolds, we design a proximal gradient method for variational mean-field games. Our comprehensive numerical experiments on various manifolds illustrate the effectiveness and flexibility of the proposed model and numerical methods.

\end{abstract}


\end{frontmatter}


\section{Introduction}
\label{sec: intro}

Mean-field games (MFG) \cite{huang2007large,huang2006large,lasry2007mean} study the behavior of a large number of rational agents in an non-cooperative game. 
It has wide applications in various fields, such as economics \cite{achdou2014partial,gomes2018mean}, engineering \cite{de2019mean,yang2017mean} as well as machine learning and reinforcement learning \cite{carmona2019model, weinan2019mean, yang2018mean,elie2020convergence}. Recently, mean field control problems have been extended into chemistry, biology , pandemic control, traffic flow models and social dynamics \cite{LLO1,LLTLO,LiLeeSO,LiLiuSO,gao2021modeling}. 
An important task in mean-field games is to study the flow of all the agents in the state space and to understand behavior of mean-field Nash equilibrium.

Conventional studies of MFG focus on choice of the state space as a Euclidean flat domain, for instance, $[0,1]^d$ with periodic boundary conditions. 
Besides research on Euclidean flat domains, there are existing work focusing MFGs on graphs \cite{GLM} or graphon state spaces ~\cite{gomes2013continuous,gueant2015existence,caines2018graphon}.
However, such spaces may not be adequate to reflect the metric structure of state spaces in many applications.
For instance, the problems of population flows or resource distributions on the Earth are actually defined on a sphere.
In machine learning, the manifold hypothesis is commonly used  \cite{cayton2005algorithms,fefferman2016testing}, since many real world data sets are actually samples from low-dimensional manifolds in a high dimensional ambient space. 
Therefore, it is quite natural and necessary to explore mean-field game/control problems on manifolds.
In this work, we would like to generalize the concepts of finite horizon mean-field games and mean-field Nash Equilibrium from Euclidean spaces to manifolds, and propose a numerical method to compute the Nash Equilibrium.

In this study, we consider a game with infinitely many indistinguishable agents on a compact and smooth manifold $\mfd$ within the time interval $[0,1]$.
At any time $t\in[0,1]$, each agent is in a certain state $\vx\in\mfd$ and the state of all agents forms a distribution $\rho(\cdot,t)\in\prob(\mfd)$. 
For each agent at $t$, given its current state $\vx$ and the anticipation of future state distribution $\rho(\cdot,s),s\in[t,1]$, the game is to optimize a control $\ctl(\vx(s),s)$ to guide its future trajectory $\vx(s),s\in(t,1]$ in order to minimize a cost $\Jfunc^{\rho}(\vx,t,\ctl)$.
Therefore the optimal control $\ctl$ depends on the state distribution $\rho$.
Although the state change of any single agent does not change $\rho(\cdot,t)$, when all the agents take the same control, the state distribution $\rho$ changes accordingly.
Thus the optimal control $\ctl$ and the state distribution $\rho$ are interdependent, and the Nash Equilibrium \cite{nash1951non,carmona2004nash}, the special pair of $(\ctl,\rho)$, is an especially interesting topic in mean-field game. 

In the conventional Euclidean setup, it has been shown that the mean-field Nash Equilibrium is the solution of a forward-backward PDE system~\cite{lasry2007mean,huang2006large,huang2007large}. We generalize this result to MFG on manifolds.
Meanwhile, for a potential mean-field game on a Euclidean domain ~\cite{lasry2007mean,cardaliaguet2015second,benamou2017variational,briani2018stable}, its optimality condition is exactly the forward-backward PDE system under mild conditions. Thus, 
the Nash Equilibrium can be obtained by searching for the stationary point of the optimization problem.
In this work, we show that the equivalence between the PDE formulation and variational formulation of mean-field game still holds on manifolds. It is worth mentioning that \cite{solomon2016entropic} studies dynamic optimal transport, a special form of potential mean-field games, on manifolds.
In this work, we consider more general forms of mean-field games on manifolds and we are interested in both the PDE and variational formulations.

There are different approaches to numerically solve mean-field games on Euclidean domains, such as finite difference methods \cite{achdou2010mean,achdou2013mean}, monotone flows \cite{almulla2017two,gomes2021numerical}, optimization algorithms \cite{benamou2015augmented,briceno2019implementation,yu2021fast} and neural networks \cite{carmona2019convergence,lin2020apac,ruthotto2020machine}. 
We refer readers to the surveys \cite{achdou2020mean,lauriere2021numerical} for more details of the numerical methods on Euclidean domains.
In our manifold setting, we focus on the variational formulation to numerically compute the Nash Equilibrium.
With the help of triangular mesh and computational geometry strategies \cite{meyer2003discrete}, we approximate the manifold, probability space and vector field space and formulate the discrete optimization problem. Once the discretization is provided, most of existing optimization-based algorithms can be adapted to solve the proposed discretization problem. In this work, we specifically use an optimization-based algorithm proposed in~\cite{yu2021fast} since it is flexible and efficient. This algorithm is adapted from the proximal gradient descent method considered in \cite{rockafellar1970convex,bauschke2011fixed,beck2009fast}. 

\textbf{Contributions:} We summarize our contributions as follows:
\begin{enumerate}
    \item[(i)] We generalize the concept of mean-field games to manifolds and derive the corresponding geometric PDE formulation of the Nash Equilibrium.
    \item[(ii)] We show the equivalence of the PDE formulation and variational formulation of mean-field games on manifolds.
    \item[(iii)] We propose a numerical method for solving the variational problem based on a proximal gradient descent method. Comprehensive experiments demonstrate the effectiveness of the proposed method.
\end{enumerate}

\textbf{Organization: } Our paper is organized as follows. 
In section~\ref{sec: mfg mf}, we derive the PDE formulation of mean-field Nash Equilibrium on manifolds. We also show that the PDE system is actually the optimality condition of an optimization problem, the potential mean-field game, on the manifold. 
We discretize the potential MFG problems in space and time domain in section~\ref{sec: disct}, 
and adapt a proximal gradient method to solve the discrete counterparts in section~\ref{sec: alg}. 
In section~\ref{sec: num res}, we provide numerical experiments that solve potential mean-field games with local or non-local interaction costs on different manifolds.

\section{Mean-field games on manifolds}
\label{sec: mfg mf}

In this section, we generalize the concepts of finite horizon mean-field games (MFGs) and their variational forms from conventional Euclidean domains to smooth and compact Riemannian manifolds. 

\subsection{Mean-field games on manifold}
\label{subsec: mfg mf}
Let's begin with some notations for convenience. We consider MFG on $(\mfd,g)$, a $d$-dimensional smooth compact manifold with a Riemannian metric $g$. As a natural extension of MFG on Euclidean domains, controls at the state $\vx\in\mfd$ are defined as elements in $\tgspx$, the tangent space of $\mfd$ at $\vx\in\mfd$. We further denote $\tgsp =\{(\vx,\vp)~|~ \vp\in\tgspx \}$ for the tangent bundle of $\mfd$; use $\Gamma(\tgsp)$ for the set of continuous vector fields on $\mfd$; and write $\prob(\mfd)$ for all probability density on $(\mfd,g)$ under the volume measurement induced by the metric $g$. 

To derive a first-order MFG system on $\mfd$, we consider a finite horizon game on the time interval $[0,1]$ with the state space $\mfd$. More specifically, we assume that there is a continuum number of agents, and each agent takes a state $\vx\in\mfd$ at any time $t\in[0,1]$. We write the state density of all the agents along $t\in[0,1]$ as $\rho\in C([0,1];\prob(\mfd))$; and assume that the impact of any single agent to $\rho$ is negligible.
Since all the agents have the same goal in a mean-field game, it is sufficient to take a representative agent as an example.
Suppose that an agent is in state $\vx$ at time $t$, the agent aims at choosing a control $\ctl\in C((t,1];\Gamma(\tgsp))$ to guide the trajectory
\begin{equation}
    \dd\vx(t) = \ctl(\vx,t)\dd t.
\label{eq: agent dynamic}
\end{equation}
in order to minimize the cost
\begin{equation}
    \Jfunc^{\rho}(\vx,t,\ctl):=\int_t^1 
    \left[ \costfuncctl(\vx(s),\ctl(\vx(s),s) + 
    \costfuncevo(\vx(s),\rho(\cdot,s)) \right] \dd s
    +\costfuncterm(\vx(1),\rho(\cdot,1)).
\label{eq: J func}
\end{equation}
Here $\costfuncctl:\tgsp\to[0,+\infty)$ is the dynamic cost, $\costfuncevo:\mfd\times \prob(\mfd)\to[0,+\infty)$ is the interaction cost, $\rho(\cdot,s)\in \prob(\mfd)$ is the density of all agents at time $s$, and $\costfuncterm:\mfd\times\prob(\mfd)\to[0,+\infty)$ is the terminal cost. 
Note that the control $\ctl$ and the state distribution $\rho$ are involved interactively. The optimal control $\ctl^{\rho}:=\argmin_{\ctl}\Jfunc^{\rho}(\vx,t,\ctl)$ generally depends on the evolution of the state distribution $\rho$.
Meanwhile, with given initial distribution $\rho(\cdot,0):=\rho_0\in\prob(\mfd)$, the distribution of agents $\rho$ is determined by the control $\ctl\in C([0,1];\Gamma(\tgsp))$ through equation \eqref{eq: agent dynamic}.
The mean-field game problem is especially interested in a special pair of them, the Nash Equilibrium, which is the same as the conventional Euclidean case \cite{nash1951non,carmona2004nash},

\begin{definition}[Nash Equilibrium]
A pair of control and state distribution $(\ctl,\rho)$ is called a Nash Equilibrium if the following two conditions hold,
\begin{enumerate}
    \item (Optimality) For any $t\in[0,1],\vx\in\mfd$,
    $J^{\rho}(\vx,t,\ctl) \leq J^{\rho}(\vx,t,\vu),~\forall~\vu\in C([0,1];\Gamma(\tgsp))$.
    \item (Consistency) $\rho(\cdot,0)=\rho_0$ where $\rho_0$ is the state distribution of all the agents at $t=0$. And $\rho(\cdot,t)$ is the state distribution of all the agents at time $t$ following the control $\ctl$.
\end{enumerate}
\end{definition}

With the definition, if $(\ctl,\rho)$ is a Nash Equilibrium of a MFG on $(\mfd, g)$, then the optimality condition ensures that $\ctl$ is the optimal control for given state distribution $\rho$, and the consistency requires that $\ctl$ lead to the state distribution $\rho$. 

According to \cite{lasry2007mean,huang2006large}, in Euclidean space, a Nash Equilibrium can be described by a PDE system, which includes a backward Hamiltonian-Jacobi-Bellman (HJB) equation by the optimality condition and a forward continuity equation by the consistency condition. In the rest, 
we will establish a similar PDE description of a Nash Equilibirum on manifolds.

Similar as the Euclidean case \cite{lasry2007mean,huang2006large}, let the value function $\vfunc$ be the cost with the optimal control,
\begin{equation}
    \vfunc^\rho(\vx,t) := \inf_{\ctl\in C([0,1];\Gamma(\tgsp))}\Jfunc^\rho (\vx,t,\ctl).
\label{eq: value func def}
\end{equation}
and $\hmtctl$ be the manifold Hamiltonian 
defined on the tangent bundle of $\mfd$ \cite{lee2017global}
\begin{equation}
    \hmtctl:\tgsp\to\bbR,\quad \hmtctl(\vx,\vq):=\sup_{\vp\in\tgspx} \left\{ -\langle\vq,\vp\rangle_{\mtc(\vx)} - \costfuncctl(\vx,\vp) \right\}.
\label{eq: control func hamiltonian def}
\end{equation}
we have the following theorem.
\begin{theorem}
If $(\ctl,\rho)$ is a Nash Equilibrium of the aforementioned mean-field game on $(\mfd,\mtc)$, then 
\begin{equation}
    \ctl(\vx,t) =\argmin_{\vp\in\tgspx} \left\{ \costfuncctl(\vx,\vp) + \langle\nabla_{\mfd}\phi(\vx,t),\vp \rangle_{\mtc(\vx)} \right\} =  -\partial_{\vq}\hmtctl(\vx,\nabla_{\mfd}\phi(\vx,t)).
\label{eq: optimal control}
\end{equation}
and $\rho,\phi$ solve
\begin{equation}
\left\{\begin{aligned}
    &-\partial_t \phi(\vx,t) +\hmtctl(\vx,\nabla_{\mfd}\phi(\vx,t)) = \costfuncevo(\vx,\rho(\cdot,t)),\\
    &\partial_t\rho(\vx,t) - \nabla_{\mfd}\cdot(\rho(\vx,t)\partial_{\vq}\hmtctl(\vx,\nabla_{\mfd}\phi(\vx,t))) = 0,\\
    &\phi(\vx,1) = \costfuncterm(\vx,\rho(\cdot,1)),\qquad\rho(\cdot,0) = \rho_0.
\end{aligned}\right.
\label{eq: mfmfg pde}
\end{equation}
\label{thm:NE_PDE}
\end{theorem}

Before proving the theorem, we give several remarks to explain notations.
\begin{remark}
\label{rem: mfd grad div}
We emphasis that the metrics $\langle\cdot,\cdot\rangle_{\mtc(\vx)}$, and operators $\nabla_{\mfd},\nabla_{\mfd}\cdot$ are based on the manifold metric $\mtc$ as a generalization of the conventional equations which only depend on the flat Euclidean metric. More details on differential geometry can be refereed in \cite{lee2013smooth}.

As an example, let $\mfd$ be a two-dimensional manifold embedded in $\bbR^3$ with an induced metric $\mtc$ on the manifold.   To be precise, consider $X:\Xi\subset\bbR^2\to\mfd\subset \bbR^3,(\xi^1,\xi^2)\mapsto X(\xi^1,\xi^2)$ as a local chart of $\mfd$, then for any $\vx=X(\xi^1,\xi^2)$ on the chart, the tangent space is $\tgspx=\spn\{\partial_{\xi^1}X,\partial_{\xi^2}X\}$. 
The matrix representation of the induced metric $g$ in the coordinate chart $X$ is provided as:
\begin{equation}
    \mtc_X(\vx):=  \left(\begin{matrix}    (\partial_{\xi^1}X)^\top \partial_{\xi^1}X & (\partial_{\xi^1}X)^\top \partial_{\xi^2}X \\ (\partial_{\xi^2}X)^\top \partial_{\xi^1}X  & (\partial_{\xi^2}X)^\top \partial_{\xi^2}X   \end{matrix}\right).
\end{equation}
Any tangent vectors $\vp,\vq\in\tgspx$ have the corresponding coordinate decomposition 
$\vp=p_X^1\partial_{\xi^1}X+p_X^2\partial_{\xi^2}X$,
$\vq=q_X^1\partial_{\xi^1}X+q_X^2\partial_{\xi^2}X$,
and the metric $g$ on each point $\vx$ is
\begin{equation}
    \mtc(\vx):\tgspx\times\tgspx\to\bbR, \quad 
    \langle\vp,\vq\rangle_{\mtc(\vx)}:=
    \left(\begin{matrix}p_X^1&p_X^2\end{matrix}\right)  \mtc_X(\vx)
    \left(\begin{matrix}q_X^1\\q_X^2\end{matrix}\right).
\end{equation}
Based on this metric, we have the following definitions of gradient and divergence operators:
\begin{eqnarray}
    \nabla_{\mfd}\phi(\vx) &=  &  
     \left(\begin{matrix}\partial_{\xi^1}\phi_X&\partial_{\xi^2}\phi_X\end{matrix}\right)  (\mtc_X(\vx))^{-1}
    \left(\begin{matrix}\partial_{\xi^1}X\\ \partial_{\xi^2}X\end{matrix}\right). \\
        \nabla_{\mfd}\cdot\ctl(\vx) &=& \frac{1}{\sqrt{\det(\mtc_X(\vx))}}\sum_{d=1}^2\frac{\partial}{\partial\xi^d}\left( \sqrt{\det(\mtc_X(\vx))}v^d_X\right).
\end{eqnarray}
where $\phi_X(\xi^1,\xi^2) = \phi(X(\xi^1,\xi^2))$ is the local representation of $\phi$ under the coordinate chart $X$, 
and the tangent vector field $\ctl\in\Gamma(\tgsp)$ has the local coordinate representation $\ctl(X(\xi^1,\xi^2))=v^1_X(\xi^1,\xi^2)\partial_{\xi^1}X+v^2_X(\xi^1,\xi^2)\partial_{\xi^2}X$. 
While the above definitions are provided in terms of a specific coordinate representation $X$, the definitions $\mtc(\vx),\nabla_{\mfd}\phi$ and $\nabla_{\mfd}\cdot\ctl$ are invariant to coordinates. 
\end{remark}

\begin{remark}
\label{rem: hamiltonian}
Following the settings and notations in previous remark, we take the quadratic dynamic cost function as an example,
\begin{equation}
\costfuncctl(\vx,\vp) := \half\|\vp\|_{\mtc(\vx)}^2=\half\langle\vp,\vp\rangle_{\mtc(\vx)}.
\end{equation}
By definition, the Hamiltonian is
\begin{equation}
\begin{aligned}
    \hmtctl(\vx,\vq)&:=\sup_{\vp\in\tgspx}\left\{-\langle\vq,\vp\rangle_{\mtc(\vx)}-\half\langle\vp,\vp\rangle_{\mtc(\vx)} \right\}\\
    &=\sup_{p^1,p^2\in\bbR}\left\{ 
    \left(\begin{matrix}-q_X^1-\half p^1 & -q_X^2-\half p^2\end{matrix}\right)    \left(\begin{matrix}    \partial_{\xi^1}X^\top\\ \partial_{\xi^2}X^\top    \end{matrix}\right)
    \left(\begin{matrix}    \partial_{\xi^1}X&\partial_{\xi^2}X    \end{matrix}\right)    \left(\begin{matrix}p^1\\p^2\end{matrix}\right)
    \right\}\\
    &=\half \left(\begin{matrix}q_X^1&q_X^2\end{matrix}\right)    \left(\begin{matrix}    \partial_{\xi^1}X^\top\\ \partial_{\xi^2}X^\top    \end{matrix}\right)
    \left(\begin{matrix}    \partial_{\xi^1}X&\partial_{\xi^2}X    \end{matrix}\right)    \left(\begin{matrix}q_X^1\\q_X^2\end{matrix}\right)\\
    &=\half\|\vq\|_{\mtc(\vx)}^2.
\end{aligned}
\end{equation}
Now we view $\tgspx=\spn\{\partial_{\xi^1}X,\partial_{\xi^2}X\}$ as a manifold and consider the nature coordinate representation $\vq=q_X^1\partial_{\xi^1}X+q_X^2\partial_{\xi^2}X$ and the induced metric $\mtc_X(\vx)$.
Then $\hmtctl(\vx,\cdot):\tgspx\to\bbR$ has the coordinate form
\begin{equation}
    \hmtctl_X(\vx,q^1,q^2) := \hmtctl(\vx,q^1\partial_{\xi^1}X + q^2\partial_{\xi^2}X)=\half \left(\begin{matrix} q^1&q^2 \end{matrix}\right)\mtc_X(\vx)
    \left(\begin{matrix} q^1\\q^2 \end{matrix}\right)
\end{equation}
and by definition of manifold gradient
\begin{equation}
    \partial_{\vq}\hmtctl(\vx,\vq) = \left(\begin{matrix}\partial_{q^1}\hmtctl_X&\partial_{q^2}\hmtctl_X\end{matrix}\right) (\mtc_X(\vx))^{-1} \left(\begin{matrix}\partial_{\xi^1}X\\ \partial_{\xi^2}X\end{matrix}\right) = \left(\begin{matrix}q_X^1&q_X^2\end{matrix}\right)\mtc_X(\vx)(\mtc_X(\vx))^{-1}
\left(\begin{matrix}\partial_{\xi^1}X\\ \partial_{\xi^2}X\end{matrix}\right) = \vq.
\end{equation}

Similarly, if we take first-order dynamic cost $\costfuncctl(\vx,\vp)=\|\vp\|_{\mtc(\vx)}$, then the corresponding Hamitonian satisfies $\hmtctl(\vx,\vq)=0$ and $\partial_{\vq}\hmtctl(\vx,\vq)=0.$
\end{remark}

\begin{remark}
\label{rem: local coord}
With manifold-metric based notations explained in remarks \ref{rem: mfd grad div} and \ref{rem: hamiltonian}, the PDE system \eqref{eq: mfmfg pde} is a generalization form to that in a Euclidean space. 
To see the difference, we state the system \eqref{eq: mfmfg pde} in a coordinate chart $X$. 
Denoting $\rho_X,\phi_X,\hmtctl_X$ and $\mtc_{\xi}$ as the local coordinate representations of $\rho, \phi, \hmtctl$ and $\mtc_X$ under $X$, respectively, we have the coordinate representation of \eqref{eq: mfmfg pde}
\begin{equation}
\left\{\begin{aligned}
    &-\partial_t \phi_X(\xi^1,\xi^2,t) +\hmtctl_X\left( X(\xi^1,\xi^2), \left(\begin{matrix}\partial_{\xi^1}\phi_X&\partial_{\xi^2}\phi_X\end{matrix}\right)\mtc_\xi^{-1}(\xi^1,\xi^2,t)
    \right) = \costfuncevo(X(\xi^1,\xi^2),\rho(\cdot,t)),\\
    &\partial_t\rho_X(\xi^1,\xi^2,t) -
    \frac{1}{\sqrt{\det(\mtc_\xi))}}  
    \sum_{i=1}^2  \frac{\partial}{\partial\xi^i}
    \left(\sum_{j=1}^2\sqrt{\det(\mtc_\xi)}\rho_X(\mtc_\xi)^{-1}_{ij}\partial_{q^{j}}\hmtctl_X\right)=0,\\
    &\phi_X(\xi^1,\xi^2,1) = \costfuncterm(X(\xi^1,\xi^2),\rho (\cdot,1)),\qquad\rho_X(\xi^1,\xi^2,0) = \rho_0(X(\xi^1,\xi^2)).
\end{aligned}\right.
\label{eq: mfmfg pde coord}
\end{equation}
with $\partial_{q^{j}}\hmtctl_X$ evaluated at 
$\left(X, \left(\begin{matrix}\partial_{\xi^1}\phi_X&\partial_{\xi^2}\phi_X\end{matrix}\right)\mtc_\xi^{-1}\right)(\xi^1,\xi^2,t) $.

It is clear to see that the above system is consistent to formula in the Euclidean case by choosing $\mfd = \bbR^2$ and $g$ as the conventional flat Euclidean metric. 
\end{remark}

Next,  we prove theorem \ref{thm:NE_PDE}.
\begin{proof}
By definition, the terminal boundary condition of $\vfunc$ is 
\begin{equation}
\phi(\vx,1) =  \costfuncterm(\vx,\rho(\cdot,1)).
\label{eq: value func term cond}
\end{equation}
For $t\in[0,1)$, by optimality of $\phi$ and dynamic programming principle, for any $h>0$
\begin{equation}
\begin{aligned}
    \phi(\vx,t) = \inf_{\ctl\in C([t,1];\Gamma(\tgsp))}
    \left\{\int_{t}^{t+h} \left[ \costfuncctl(\vx(s),\ctl(\vx(s),s)) + 
    \costfuncevo(\vx(s),\rho(\cdot,s)) \right] \dd s 
    +\phi(\vx(t+h),t+h)\right\}.
\end{aligned}
\label{eq: value func dynamic programming}
\end{equation} 
where $\vx(t+h) = \vx(t)+\int_{t}^{t+h}\ctl(\vx(s),s)\dd s.$
Assume that $\phi$ is $C^2$ in $\vx$ and $C^1$ in $t$. Then by Leibniz integral rule, 
\begin{equation}
\begin{aligned}
    \phi(\vx(t+h),t+h) = \phi(\vx,t) 
    + \int_{t}^{t+h} \left[ \partial_t\phi(\vx(s),s) + \left\langle\nabla_{\mfd}\phi(\vx(s),s),\ctl(\vx(s),s)\right\rangle_{\mtc(\vx)} \right] \dd s,
\end{aligned}
\label{eq: value func leibniz rule}
\end{equation}
Combining \eqref{eq: value func dynamic programming} and \eqref{eq: value func leibniz rule}, we have
\begin{equation}
\begin{aligned}
    &\int_{t}^{t+h} \left[ \partial_t\phi(\vx(s),s) + \costfuncevo(\vx(s),\rho(\cdot,s)) \right] \dd s \\
    &+ \inf_{\ctl\in C([t,t+h];\tgsp)} \left\{ \int_t^{t+h}\left[ \costfuncctl(\vx(s),\ctl(\vx(s),s)) + \langle\nabla_{\mfd}\phi(\vx(s),s),\ctl(\vx(s),s)\rangle_{\mtc(\vx)} \right] \dd s \right\} = 0.
\end{aligned}
\label{eq: value func hjb integral form}
\end{equation}
Dividing both sides by $h$ and letting $h\to 0^+$ shows that $\phi^\rho$ satisfies the HJB equation \eqref{eq: value func hjb no bd} on $\mfd$.
\begin{equation}
    -\partial_t \phi(\vx,t) - \inf_{\vp\in\tgspx} \left\{ \costfuncctl(\vx,\vp) + \langle\nabla_{\mfd}\phi(\vx,t),\vp\rangle_{\mtc(\vx)} \right\} = \costfuncevo(\vx,\rho(\cdot,t)).
\label{eq: value func hjb no bd}
\end{equation}
Plugging in the definition of manifold Hamiltonian
\begin{equation}
    \hmtctl:\tgsp\to\bbR,\quad \hmtctl(\vx,\vq)=\sup_{\vp\in\tgspx} \left\{ - \costfuncctl(\vx,\vp) -\langle\vq,\vp\rangle_{\mtc(\vx)}  \right\}.
\end{equation}
we show that $\phi$ satisfies the HJB equation
\begin{equation}
\left\{\begin{aligned}
    &-\partial_t \phi(\vx,t) +\hmtctl(\vx,\nabla_{\mfd}\phi(\vx,t)) = \costfuncevo(\vx,\rho(\cdot,t)),\\
    &\phi(\vx,1) = \costfuncterm(\vx,\rho(\cdot,1)).
\end{aligned}\right.
\label{eq: hjb}
\end{equation}
And by properties of convex conjugate, we obtain the optimal control
\begin{equation}
    \ctl(\vx,t) :=\argmin_{\vp\in\tgspx} \left\{ \costfuncctl(\vx,\vp) + \langle\nabla_{\mfd}\phi(\vx,t),\vp\rangle_{\mtc(\vx)} \right\} =  -\partial_{\vq}\hmtctl(\vx,\nabla_{\mfd}\phi(\vx,t)).
\end{equation}

On the other hand, by consistency condition of a Nash Equilibrium, with initial density $\rho_0$, $\rho$ satisfies the continuity equation driven by $\ctl(\vx,t)$, 
\begin{equation}\left\{
\begin{aligned}
    &\partial_t\rho(\vx,t) + \nabla_{\mfd}\cdot(\rho(\vx,t)\ctl(\vx,t)) = 0,\\
    &\rho(\cdot,0) = \rho_0.
\end{aligned}\right.
\end{equation}
And $\ctl$ being the optimal control $\ctl(\vx,t)=-\partial_{\vq}\hmtctl(\vx,\nabla_{\mfd}\phi(\vx,t))$ implies that $\rho$ satisfies
\begin{equation}\left\{
\begin{aligned}
    &\partial_t\rho(\vx,t) - \nabla_{\mfd}\cdot(\rho(\vx,t)\partial_{\vq}\hmtctl(\vx,\nabla_{\mfd}\phi(\vx,t))) = 0,\\
    &\rho(\cdot,0) = \rho_0.
\end{aligned}\right.
\label{eq: cteqn}
\end{equation}

To summarize, the solution $(\phi,\rho)$ to the following PDE system gives us a Nash Equilibrium $(\ctl,\rho)$ with $\ctl=-\partial_{\vq}\hmtctl(\vx,\nabla_{\mfd}\phi(\vx,t))$,
\begin{equation}
\left\{\begin{aligned}
    &-\partial_t \phi(\vx,t) +\hmtctl(\vx,\nabla_{\mfd}\phi(\vx,t)) = \costfuncevo(\vx,\rho(\cdot,t)),\\
    &\partial_t\rho(\vx,t) - \nabla_{\mfd}\cdot(\rho(\vx,t)\partial_{\vq}\hmtctl(\vx,\nabla_{\mfd}\phi(\vx,t))) = 0,\\
    &\phi(\vx,1) = \costfuncterm(\vx,\rho(\cdot,1)),\qquad\rho(\cdot,0) = \rho_0.
\end{aligned}\right.
\end{equation}
\end{proof}

At the end of this part, we present some common examples. 
\begin{example}[Local mean-field games]
\label{eg: local cost}
When the interaction cost and terminal cost functions take the local form. i.e. the cost at $\vx$ only depends on the density at $\vx$. 
the corresponding mean-field game is called a local mean-field game.

We list some common choices of $\costfuncevo$ and $\costfuncterm$ here.
\begin{itemize}
    \item $\costfuncevo(\vx,\rho(\cdot,t))=B(\vx)$ with $B:\mfd\to\bbR$. 
    This interaction function gives a preference of states. The agents tend to stay at $\vx$ where the cost $B(\vx)$ is low. 
    \item $\costfuncevo(\vx,\rho(\cdot,t))=\log(\rho(\vx,t))+1$ and
    $\costfuncevo(\vx,\rho(\cdot,t))=(\rho(\vx,t))^p,p>0$.
    These interaction functions discourage the aggregation of densities.
    \item $\costfuncterm(\vx,\rho(\cdot,1)) = (\rho(\vx,1)-\rho_1(\vx))^2$ and  $\costfuncterm(\vx,\rho(\cdot,1)) = \log(\rho(\vx,1))-\log(\rho_1(\vx))+1$ with a given $\rho_1$.
    These terminal functions encourage $\rho(\cdot,1)$ to approach to the desired terminal density $\rho_1$.
\end{itemize}
\end{example}

\begin{example}[Non-local mean-field games]
\label{eg: nonlocal cost}
The interaction cost function $\costfuncevo$ or terminal cost function $\costfuncterm$ can also take non-local forms. Take $\costfuncevo$ as an example. 
If $\kernel:\mfd\times\mfd\to\bbR$ is a convolutional kernel, and
\begin{equation}
    \costfuncevo(\vx,\rho(\cdot,t)):=\int_{\mfd}\kernel(\vx,\vy)\rho(\vy,t)\mathrm{d}_\mfd \vy,
\end{equation}
then the mean-field game is non-local.
Symmetric kernel functions $\kernel$ with $\kernel(\vx,\vy)=\kernel(\vy,\vx)$ are of special interest to us. When $\kernel$ is symmetric, the PDE system is the optimality condition of a variational problem \cite{nurbekyan2019fourier,liu2021computational}. We provide detailed discussions in the following section. \end{example}

\subsection{Potential mean-field games on manifold}
\label{subsec: vmfg mf}

According to \cite{lasry2007mean,cardaliaguet2015second,benamou2017variational,briani2018stable}, when the state space is Euclidean, with mild conditions, the local minimizer of an optimization problem and the corresponding dual variable is a weak solution to the MFG PDE system. 
In this part, we establish the parallel results on manifolds. 
We formulate the potential MFG on manifold and show that the necessary optimality condition of this variational problem is exactly the PDE system \eqref{eq: mfmfg pde} under certain conditions.

\begin{theorem}
\label{thm: potential=pde}
Assume that $\costfuncctl(\vx,\vp)$ is convex in $\vp\in\tgspx$ at any $\vx\in\mfd$, and there exist $\Costevo: \prob(\mfd)\to[0,+\infty)$, $\Costterm: \prob(\mfd)\to[0,+\infty)$ such that $\frac{\delta\Costevo(\rho)}{\delta\rho}(\vx) = \costfuncevo(\vx,\rho), \frac{\delta\Costterm(\rho)}{\delta\rho}(\vx) = \costfuncterm(\vx,\rho)$.
Consider the optimization problem,
\begin{equation}
\begin{aligned}
    \inf_{\rho,\vm}\quad &\Cost(\rho,\vm):=\int_0^1\int_{\mfd} \rho(\vx,t)\costfuncctl\left(\vx,\frac{\vm(\vx,t)}{\rho(\vx,t)}\right)\ddx\dd t + \int_0^1\Costevo(\rho(\cdot,t))\dd t + \Costterm(\rho(\cdot,1))\\
    \text{ subject to } &\quad \partial_t\rho + \nabla_{\mfd}\cdot\vm = 0, \rho(\cdot,0) = \rho_0,\\
    &\quad \rho\in C([0,1];\prob(\mfd)),\vm\in C([0,1];\Gamma(\tgsp)).
\end{aligned}
\label{eq: mfmfg opt}
\end{equation}
When $\rho(\vx,t)=0$, we take the conventional definition of $\costfuncctl$
\begin{equation}
    \costfuncctl\left(\vx,\frac{\vm(\vx,t)}{\rho(\vx,t)}\right)=
    \begin{cases}
    0,&\text{if}\quad \vm(\vx,t)=\Zero,\\
    +\infty,&\text{if} \quad \vm(\vx,t)\neq\Zero.
    \end{cases}
\end{equation}
The following statements hold
\begin{enumerate}
    \item If $(\rho,\vm)$ is a local minimizer of \eqref{eq: mfmfg opt}, then there exists $\phi$ such that $\vm=-\rho\partial_{\vq}\hmtctl(\cdot,\nabla_{\mfd}\phi(\vx,t))$ and $(\rho,\phi)$ solves
\begin{equation}
\left\{\begin{aligned}
    &-\partial_t \phi(\vx,t) +\hmtctl(\vx,\nabla_{\mfd}\phi(\vx,t)) \leq \costfuncevo(\vx,\rho(\cdot,t)),\\
    &\partial_t\rho(\vx,t) - \nabla_{\mfd}\cdot(\rho(\vx,t)\partial_{\vq}\hmtctl(\vx,\nabla_{\mfd}\phi(\vx,t))) = 0,\\
    &\phi(\vx,1) \leq \costfuncterm(\vx,\rho(\cdot,1)),\qquad\rho(\cdot,0) = \rho_0.
\end{aligned}\right.
\label{eq: mfmfg kkt}
\end{equation}
In addition, if $\rho>0$, then $(\rho,\phi)$ solves the PDE system \eqref{eq: mfmfg pde}.

    \item If $\Cost$ is pseudo-convex in $(\rho,\vm)$, and $(\phi,\rho)$ is a solution to the MFG PDE system \eqref{eq: mfmfg pde}, then $\vm=-\rho\partial_{\vq}\hmtctl(\cdot,\nabla_{\mfd}\phi(\vx,t))$, and $(\rho,\vm)$ is the minimizer of \eqref{eq: mfmfg opt}.
\end{enumerate}

\end{theorem}

\begin{proof}

We first derive the KKT system of \eqref{eq: mfmfg opt} based on the theory of constrained optimization ~\cite{kuhn1951nonlinear}.
We denote $\phi\in C([0,1]\times\mfd)$ as the Lagrangian multiplier for the continuity equation, and then the Lagrangian of \eqref{eq: mfmfg opt} is,
\begin{align}
    \Lag(\rho,\vm,\phi):=&\int_0^1\int_\mfd  \rho(\vx,t)\costfuncctl\left(\vx,\frac{\vm(\vx,t)}{\rho(\vx,t)}\right)\ddx\dd t + \int_0^1\Costevo(\vx,\rho(\cdot,t))\dd t
    +  \Costterm(\vx,\rho(\cdot,1)) \nonumber\\
    &- \int_0^1\int_\mfd \phi(\vx,t)\left( \partial_t\rho + \nabla_{\mfd}\cdot\vm \right)(\vx,t)\ddx\dd t\\
    =&\int_0^1\int_\mfd  \rho(\vx,t)\costfuncctl\left(\vx,\frac{\vm(\vx,t)}{\rho(\vx,t)}\right) \ddx\dd t + \int_0^1\Costevo(\vx,\rho(\cdot,t))\dd t   \nonumber\\
    &+ \int_0^1\int_\mfd \left[\rho(\vx,t) \partial_t\phi(\vx,t) + \langle\vm(\vx,t), \nabla_{\mfd}\phi(\vx,t)\rangle_\mtc \right] \ddx\dd t \nonumber\\
    &+\Costterm(\vx,\rho(\cdot,1)) + \int_\mfd  \left[-\phi(\vx,1)\rho(\vx,1)+\phi(\vx,0)\rho_0(\vx) \right] \ddx,
\end{align}
Since $\rho\geq0,\rho(\cdot,0)=\rho_0$, the KKT system of \eqref{eq: mfmfg opt} is 
\begin{equation}\left\{
\begin{aligned}
    &\delta_{\rho}\Lag(\rho,\vm,\phi) \geq 0,\quad \rho\delta_{\rho}\Lag(\rho,\vm,\phi)=0,\\
    &\delta_{\vm}\Lag(\rho,\vm,\phi) = 0,\\
    &\partial_t\rho(\vx,t) + \nabla_{\mfd}\cdot\vm(\vx,t)=0,\rho(\cdot,0)=\rho_0.\\
\end{aligned}\right.
\label{eq: mfmfg kkt unsimplified}
\end{equation}
Among the system \eqref{eq: mfmfg kkt unsimplified}, $\delta_{\vm}\Lag(\rho,\vm,\phi) = 0$ yields
\begin{equation}
    \partial_\vp\costfuncctl\left(\vx,\frac{\vm(\vx,t)}{\rho(\vx,t)}\right)+\nabla_{\mfd}\phi(\vx,t) = 0.
\end{equation}
and consequently $\vm=-\rho\partial_\vq\hmtctl(\cdot,\nabla_\mfd\phi)$ by convexity of $\costfuncctl$.
Plugging in and simplifying $\delta_{\rho}\Lag(\rho,\vm,\phi) \geq 0$ then gives
\begin{equation}\left\{
\begin{aligned}
    &-\partial_t\phi(\vx,t)+\hmtctl(\vx,\nabla_{\mfd}\phi(\vx,t))\leq    \costfuncevo(\vx,\rho(\cdot,t)) ,\\
    &\phi(\vx,1) \leq \costfuncterm(\vx,\rho(\cdot,1)),\\
\end{aligned}\right.
\end{equation}
and the equality hold when $\rho(\vx,t)>0$.
Combining above, we see \eqref{eq: mfmfg kkt} is exactly the KKT system \eqref{eq: mfmfg kkt unsimplified}.

According to optimization theory~\cite{kuhn1951nonlinear}, because the constraints of \eqref{eq: mfmfg opt} are linear in $(\rho,\vm)$, the KKT conditions are necessary for the local minimizer and thus the first statement holds.
In addition, when $\Cost$ is pseudo-convex, the KKT conditions are sufficient for the minimizer~\cite{mangasarian1975pseudo}.  Since the PDE system \eqref{eq: mfmfg pde} implies the KKT system of \eqref{eq: mfmfg opt}, the second statement holds.
\end{proof}

With the above theorem, the forward-backward system \eqref{eq: mfmfg pde} can be solved by searching for the local minimizer of variational problem \eqref{eq: mfmfg opt}. In this study, we majorly focus on the variational problem \eqref{eq: mfmfg opt}.

In the rest part of this section, we present some examples of potential mean-field games as well as their corresponding PDE systems. 

\begin{example}[Quadratic dynamic cost with local interaction]
\label{eg: vmfg local}
Let $\costfuncctl(\vx,\vp) = \half\|\vp\|_{\mtc(\vx)}^2$ and the local interaction and terminal costs
\begin{equation}
\begin{aligned}
    &\Costevo(\rho(\cdot,t)) = \int_{\mfd}\rho(\vx,t)\log(\rho(\vx,t))\ddx,\\
    &\Costterm(\rho(\cdot,1)) = \int_{\mfd}\rho(\vx,1)\log\left(\frac{\rho(\vx,1)}{\rho_1(\vx)}\right)\ddx.
\end{aligned}
\end{equation}
where $\rho_1(\vx)$ is a given density. 
With these choices of costs, $\costfuncctl(\vx,\cdot)$ is convex in $\vp\in\tgspx$ and the objective function $\Cost$ is pseudo-convex in $(\rho,\vm)$.
According to theorem \ref{thm: potential=pde}, searching for the optimizer is equivalent to solving a mean-field game PDE system.
To be precise, if the optimizer $\rho>0$, then the KKT system of this potential game is
\begin{equation}
\left\{\begin{aligned}
    &-\partial_t \phi(\vx,t) +\half\left\|\nabla_{\mfd}\phi(\vx,t)\right\|_{\mtc(\vx)}^2 = \log(\rho(\vx,t))+1,\\
    &\partial_t\rho(\vx,t) - \nabla_{\mfd}\cdot(\rho(\vx,t)\nabla_{\mfd}\phi(\vx,t)) = 0,\\
    &\phi(\vx,1) = \log\left(\frac{\rho(\vx,1)}{\rho_1(\vx)}\right)+1,\qquad\rho(\cdot,0) = \rho_0.
\end{aligned}\right.
\end{equation}
and $\ctl=-\nabla_{\mfd}\phi(\vx,t)$.
It is easy to check that this system is the PDE formulation of the mean-field game with
\begin{equation}
\begin{aligned}
    &\costfuncevo(\vx,\rho(\cdot,t)) = \frac{\delta\Costevo(\rho)}{\delta\rho}(\vx) = \log(\rho(\vx,t))+1,\\
    &\costfuncterm(\vx,\rho(\cdot,1)) = \frac{\delta\Costterm(\rho)}{\delta\rho}(\vx) = \log\left(\frac{\rho(\vx,1)}{\rho_1(\vx)}\right)+1.
\end{aligned}
\end{equation}

\end{example}

\begin{example}[Quadratic dynamic cost with non-local interaction cost]
\label{eg: vmfg nonlocal}
Let $\costfuncctl(\vx,\vp) = \half\|\vp\|_{\mtc(\vx)}^2$ and the local terminal costs be
\begin{equation}
    \Costterm(\rho(\cdot,1)) = \int_{\mfd}\half(\rho(\vx,1)-\rho_1(\vx))^2\ddx.
\end{equation}
with a given density $\rho_1$. 
We consider a non-local interaction cost
\begin{equation}
    \Costevo(\rho(\cdot,t))=\half\int_{\mfd\times\mfd}\kernel(\vx,\vy)\rho(\vx,t)\rho(\vy,t) \ddx \ddy,
\end{equation}
where $\kernel(\vx,\vy) = \mu\exp\left(-\frac{1}{\sigma^2}d^2_{\mtc}(\vx,\vy)\right)$ is a Gaussian kernel based on the geodesic distance $d_{\mtc}$ of $(\mfd,\mtc)$.
The KKT system of this variational problem is
\begin{equation}
\left\{\begin{aligned}
    &-\partial_t \phi(\vx,t) +\half\left\|\nabla_{\mfd}\phi(\vx,t)\right\|_{\mtc(\vx)}^2 = \int_{\mfd} \kernel(\vx,\vy)\rho(\vy,t)\ddy,\\
    &\partial_t\rho(\vx,t) - \nabla_{\mfd}\cdot(\rho(\vx,t)\nabla_{\mfd}\phi(\vx,t)) = 0,\\
    &\phi(\vx,1) = \rho(\vx,1)-\rho_1(\vx),\qquad\rho(\cdot,0) = \rho_0.
\end{aligned}\right.
\end{equation}
and $\ctl=-\nabla_{\mfd}\phi(\vx,t)$.
This system is the PDE formulation of the mean-field game with
\begin{equation}
\begin{aligned}
    &\costfuncevo(\vx,\rho(\cdot,t)) = \frac{\delta\Costevo(\rho)}{\delta\rho}(\vx) = \int_{\mfd}\kernel(\vx,\vy)\rho(\vy,t)\ddy,\\
    &\costfuncterm(\vx,\rho(\cdot,1)) = \frac{\delta\Costterm(\rho)}{\delta\rho}(\vx) = \rho(\vx,1)-\rho_1(\vx).
\end{aligned}
\end{equation}
Here $\frac{\delta\Costevo(\rho)}{\delta\rho}(\vx) = \int_{\mfd}\kernel(\vx,\vy)\rho(\vy,t)\ddy$ holds because $K$ is symmetric. For general non-symmetric kernels, the corresponding interaction cost function can be written as $\costfuncevo(\vx,\rho(\cdot,t))=\int_{\mfd}\left(\half \kernel(\vx,\vy)+\half \kernel(\vy,\vx)\right)\rho(\vy,t)\ddy$.
\end{example}

\section{Discretization on Manifolds}
\label{sec: disct}

The optimization problem \eqref{eq: mfmfg opt} is defined in an infinite-dimension space and in general is solved approximately using appropriate discretization. Different from conventional MFG problems in Euclidean space, we need to approximate the ground manifold as well as functions and vector fields on the manifold. In this section, we focus on two-dimensional manifolds and discuss the discrete counterpart of \eqref{eq: mfmfg opt}.
We first approximate the manifold with a triangular mesh. This leads a semi-dicrete version of \eqref{eq: mfmfg opt} and its associated KKT system. After that, we derive a fully discretized version for our numerical implementation by equally splitting the time interval. 

\subsection{Space discretization}
\label{subsec: vmfg trg}

\begin{figure}[h]
\centering
\includegraphics[scale=1]{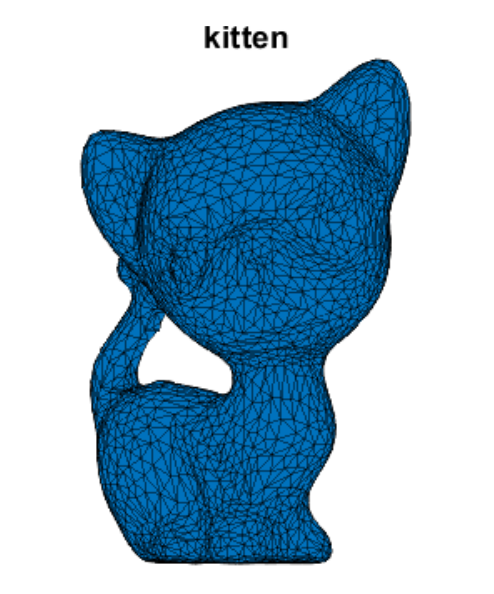}
\includegraphics[scale=1]{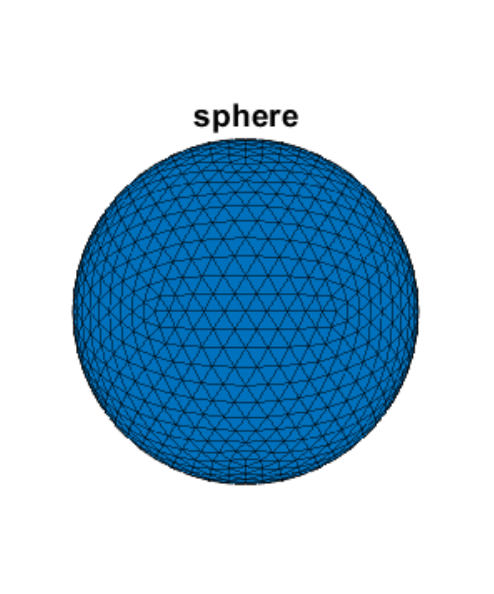}
\includegraphics[scale=1]{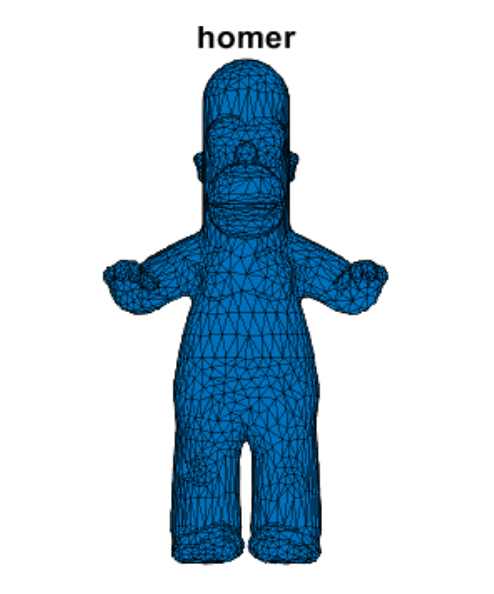}
\caption{Triangular mesh approximation of some manifolds.}
\label{fig: trimesh illu}
\end{figure}

We follow a conventional approach to approximate a two-dimensional manifold $\mfd\subset\bbR^3$ by a triangular mesh $\disctmfd:=\{\vtc,\trg\}$~\cite{meyer2003discrete,lai2011framework}. Here $\vtc=\{\vtci\in\bbR^3\}_{i=1}^{\numvtc}$ represents the set of vertices,  $\trg=\{\trgj\}_{j=1}^{\numtrg}$ is the set of triangles where each triangle has three vertices in $\vtc$. Figure \ref{fig: trimesh illu} shows several triangular meshes used in our numerical experiments later. For convenience, we also abuse our notation $\disctmfd$ for the piecewise linear approximation of $\mfd$ obtained from the the given triangular mesh. 

For a real-valued function $\psi:\mfd\to\bbR$, we approximate it with a piece-wise linear function $\Psi:\disctmfd\to\bbR$, where on vertices $\Psi(\vtci):=\psi(\vtci)$ and on each triangle $\Psi(\vx)$ is linear. 
In this way, any piece-wise linear function on $\disctmfd$ is fully represented by its values on vertices. With a slight abuse of notations, we denote $\Psi:=\left(\begin{matrix} \Psi(\vtc_1)\\ \vdots \\ \Psi(\vtc_{\numvtc}) \end{matrix}\right)$ as a vector in $\bbR^{\numvtc}$. 
By piece-wise linearity, the gradient of $\Psi$ is a piece-wise constant vector field.
For consistency, we let $\discttgsp:=\sqcup_{j=1}^{\numtrg}\spn{\trgj}$ mimic the tangent bundle,
$\Gamma(\discttgsp):=\left\{U:\trg\to\discttgsp,U(\trgj) = \left(\begin{matrix}U^1(\trgj)\\U^2(\trgj)\\U^3(\trgj)\end{matrix}\right)\in\spn{\trgj}\right\}$ denote the set of piece-wise constant vector field.
Similar to the function discretization, we use the matrix $U=\left(\begin{matrix}U^1,U^2,U^3\end{matrix}\right)=\left(\begin{matrix} (U(\trg_1))^\top\\ \vdots \\ (U(\vtc_{\numtrg}))^\top \end{matrix}\right)\in\bbR^{\numtrg\times3}$ to fully describe the vector field.

\begin{figure}[h]
\centering
\includegraphics[width=7cm]{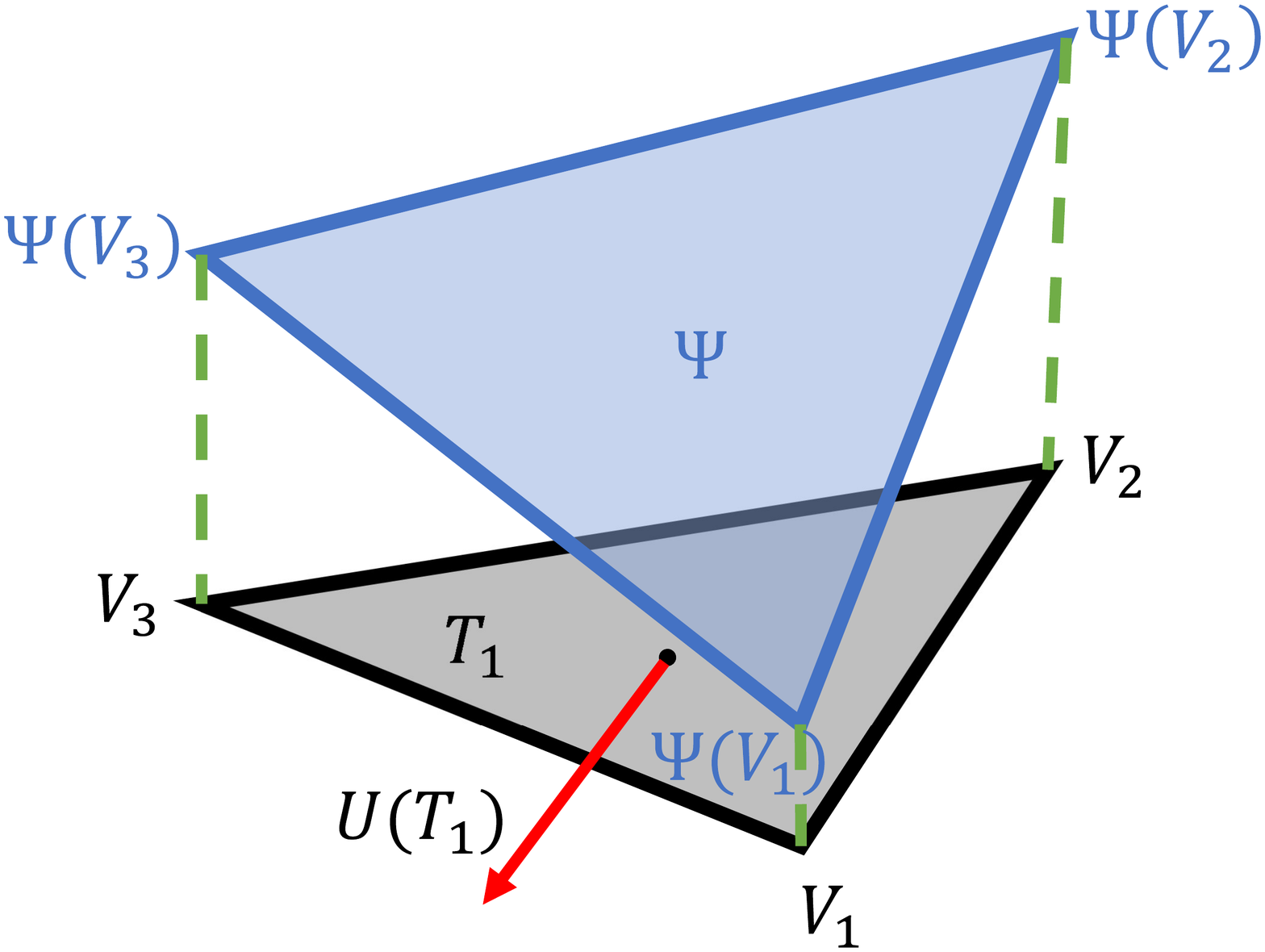}
\includegraphics[width=7cm]{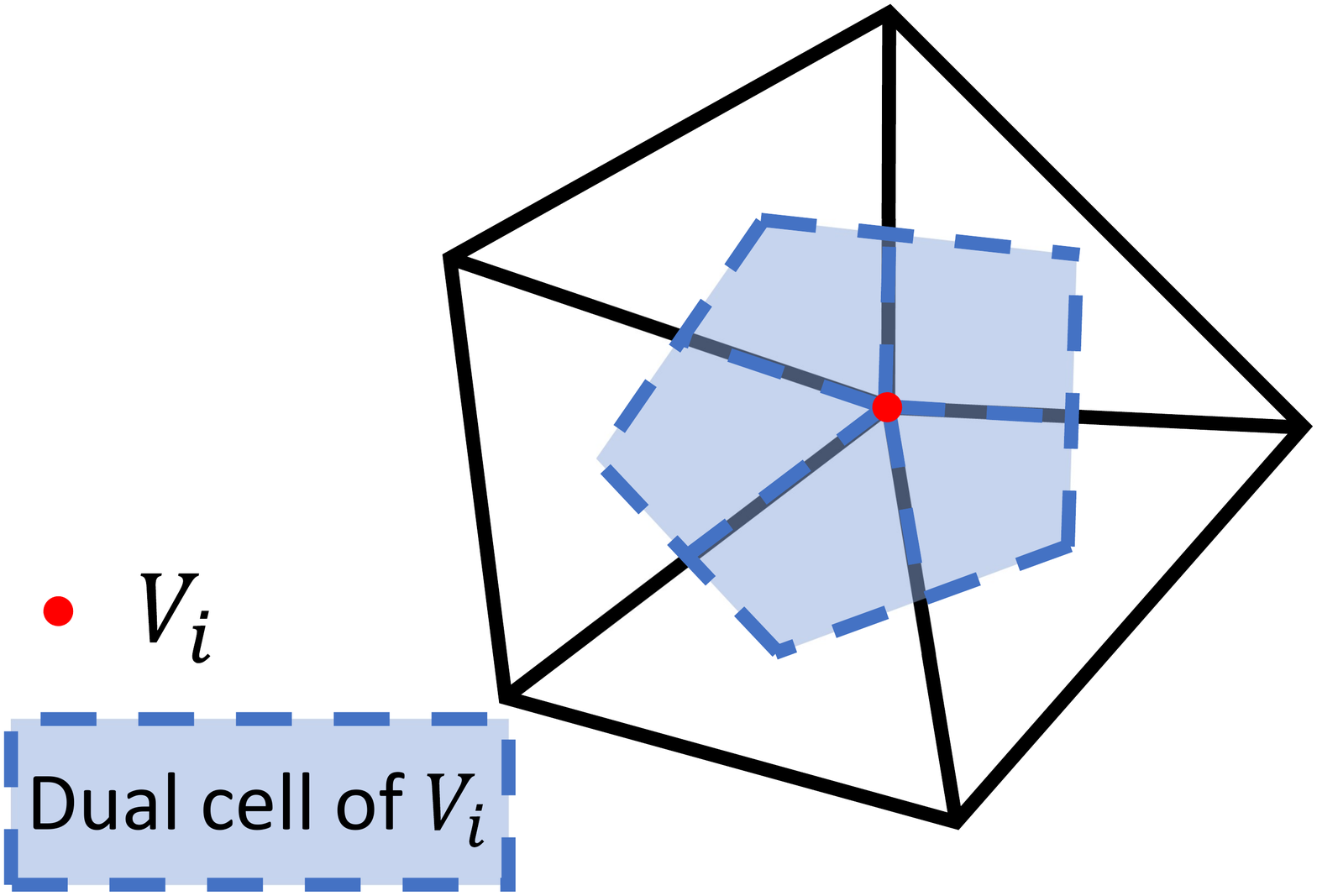}
\caption{Illustration of definition of areas, functions and vector field on triangular meshes. (Left: gradient operator on one triangular face, right: dual cell of $\vtci$.)}
\label{fig: disct illu}
\end{figure}

Given any $\Psi:\disctmfd\to\bbR$, its gradient $\grad\Psi\in\Gamma(\discttgsp)$ can be written as  $\grad\Psi=\left(\begin{matrix}\disctmtc^1\Psi&\disctmtc^2\Psi&\disctmtc^3\Psi\end{matrix}\right)$,  where $\disctmtc^1,\disctmtc^2,\disctmtc^3\in\bbR^{\numtrg\times\numvtc}$ provides a discretization of $\nabla_\mfd$.
To see this, take a triangle  $\trg_1$ with vertices $\vtc_1,\vtc_2,\vtc_3$ as an example (see the left image in Figure \ref{fig: disct illu}).
We first parameterise $\trg_1$ by
\begin{equation}
    \vtc(\xi^1,\xi^2):=\xi^1(\vtc_2-\vtc_1)+\xi^2(\vtc_3-\vtc_1),\quad 0\leq\xi^1,\xi^2\leq1,\xi^1+\xi^2=1.
\end{equation}
Thus the induced metric on $\trg_1$ is the constant matrix
\begin{equation}
    \mtc= \left(\begin{matrix} 
    \langle \vtc_2-\vtc_1,\vtc_2-\vtc_1\rangle & \langle \vtc_2-\vtc_1,\vtc_3-\vtc_1\rangle\\
    \langle \vtc_3-\vtc_1,\vtc_2-\vtc_1\rangle & \langle \vtc_3-\vtc_1,\vtc_3-\vtc_1\rangle 
    \end{matrix}\right)
\end{equation}
Because $\Psi$ restricted on $\trg_1$ is linear, we have
\begin{equation}
    \Psi(\vtc(\xi^1,\xi^2)):=\xi^1(\Psi(\vtc_2)-\Psi(\vtc_1))+\xi^2(\Psi(\vtc_3)-\Psi(\vtc_1)).
\end{equation}
The definition of gradient gives us
\begin{equation}
\begin{aligned}
    \left(\begin{matrix} 
    \langle \vtc_2-\vtc_1,(\grad\Psi)(\trg_1) \rangle\\
    \langle \vtc_3-\vtc_1,(\grad\Psi)(\trg_1) \rangle 
    \end{matrix}\right)
    &= \left(\begin{matrix} \Psi(\vtc_2)-\Psi(\vtc_1)\\ \Psi(\vtc_3)-\Psi(\vtc_1) \end{matrix}\right)
    = \left(\begin{matrix} -1 & 1 & 0\\ -1 & 0 & 1 \end{matrix}\right)
    \left(\begin{matrix} \Psi(\vtc_1)\\\Psi(\vtc_2)\\\Psi(\vtc_3)\end{matrix}\right).
\end{aligned}
\label{eq: gd op 1}
\end{equation}
Since $(\grad\Psi)(\trg_1)\in\spn{\trg_1}$, it is clear that the gradient of $\Psi$ on $\trg_1$ has the decomposition 
$ (\grad\Psi)(\trg_1) = 
\left(\begin{matrix} \vtc_2-\vtc_1,&\vtc_3-\vtc_1 \end{matrix}\right)
\left(\begin{matrix} \mu^1\\ \mu^2 \end{matrix}\right)$
and therefore we have
\begin{equation}
    \left(\begin{matrix} 
    \langle \vtc_2-\vtc_1,(\grad\Psi)(\trg_1) \rangle\\
    \langle \vtc_3-\vtc_1,(\grad\Psi)(\trg_1) \rangle 
    \end{matrix}\right)
    = \mtc\left(\begin{matrix} \mu^1\\ \mu^2 \end{matrix}\right)
\label{eq: gd op 2}
\end{equation}
Because the triangle is non-degenerative, \eqref{eq: gd op 1} and \eqref{eq: gd op 2} together solves $\mu^1,\mu^2$
\begin{equation}
\begin{aligned}
    \left(\begin{matrix} \mu^1\\ \mu^2 \end{matrix}\right)
    = \mtc^{-1}
    \left(\begin{matrix} -1 & 1 & 0\\ -1 & 0 & 1 \end{matrix}\right)
    \left(\begin{matrix} \Psi(\vtc_1)\\\Psi(\vtc_2)\\\Psi(\vtc_3)\end{matrix}\right).
\end{aligned}
\end{equation}
And this implies 
\begin{equation}
\begin{aligned}
    (\grad\Psi)(\trg_1) = 
    \left(\begin{matrix} \vtc_2-\vtc_1,&\vtc_3-\vtc_1 \end{matrix}\right)
    \mtc^{-1}
    \left(\begin{matrix} -1 & 1 & 0\\ -1 & 0 & 1 \end{matrix}\right)
    \left(\begin{matrix} \Psi(\vtc_1)\\\Psi(\vtc_2)\\\Psi(\vtc_3)\end{matrix}\right).
\end{aligned}
\end{equation}
Assigning $\disctmtc^d(\trg_1,\vtci) = 0$ for $\vtci\in \vtc\backslash\{\vtc_1,\vtc_2,\vtc_3\},d=1,2,3$ and 
\begin{equation}
\begin{aligned}
\left(\begin{matrix}
\disctmtc^1(\trg_1,\vtc_1) & \disctmtc^1(\trg_1,\vtc_2) & \disctmtc^1(\trg_1,\vtc_3) \\
\disctmtc^2(\trg_1,\vtc_1) & \disctmtc^2(\trg_1,\vtc_2) & \disctmtc^2(\trg_1,\vtc_3) \\
\disctmtc^3(\trg_1,\vtc_1) & \disctmtc^3(\trg_1,\vtc_2) & \disctmtc^3(\trg_1,\vtc_3) 
\end{matrix}\right)
= 
\left(\begin{matrix} \vtc_2-\vtc_1,&\vtc_3-\vtc_1 \end{matrix}\right)
\mtc^{-1}
\left(\begin{matrix} -1 & 1 & 0\\ -1 & 0 & 1 \end{matrix}\right)
\end{aligned}
\end{equation}
assures $(\grad\Psi)(\trg_1) = \left(\begin{matrix}(\disctmtc^1\Psi)(\trg_1)&(\disctmtc^2\Psi)(\trg_1)&(\disctmtc^3\Psi)(\trg_1) \end{matrix}\right)$. 
Following the same approach to define $\disctmtc^d(\trgj,\vtci)$ on $\trgj\in\trg$, we have $\grad\Psi=\left(\begin{matrix}\disctmtc^1\Psi&\disctmtc^2\Psi&\disctmtc^3\Psi\end{matrix}\right)$.

Next, we define discretization of the divergence operator based on its adjoint relation to the gradient operator. Consider the following discretization  of surface area and inner product.
Let $\areatrgj$ be the area of triangle $\trgj$ and $\areavtci:=\third\sumtrg\areatrgj$ be the area of the barycentric dual cell of $\vtci$ (Figure \ref{fig: disct illu} right),
and denote $\areavtc:=\diag(\area_{\vtc_1},\cdots,\area_{\vtc_{\numvtc}})\in\bbR^{\numvtc\times \numvtc}$ and $\areatrg:=\diag(\area_{\trg_1},\cdots,\area_{\trg_{\numtrg}})\in\bbR^{\numtrg\times \numtrg}$ be the mass matrices of vertices and of triangles.
We then define the inner products of vector fields as $\langle U_1,U_2 \rangle_{\trg}:=\trace(U_1^\top\areatrg U_2)=\sum_{d=1}^3 (U_1^d)^\top\areatrg U_2^d$ and of functions as $\langle\Psi_1,\Psi_2\rangle_{\vtc}:=\Psi_1^\top\areavtc\Psi_2$.
To preserve the adjoint relation between negative gradient and divergence under above inner products, i.e.  $\langle -\grad\Psi,U\rangle_\trg = \langle \Psi,\divg U\rangle_\vtc$, we assign $\divg U:=-\sum_{d=1}^3\areavtc^{-1}(\disctmtc^d)^\top\areatrg U^d\in\bbR^{\numvtc}$.

With the notations of area, we write the set of probability density functions on $\disctmfd$ as $\prob(\disctmfd):=\{\Rho\in\bbR_+^{\numvtc}:\areavtc\Rho = 1\}$. 
If the the initial density is given as $\Rho_0\in\prob(\disctmfd)$, the semi-discrete constraint is then
\begin{equation}
\begin{aligned}
    & \frac{\dd}{\dd t}\Rho(\vtci,t) + \left(\divg M\right)(\vtci,t) = 0,\Rho(\cdot,0)=\Rho_0,\\
    & \Rho\in C([0,1];\prob(\vtc)), M\in C([0,1];\Gamma(\discttgsp)).
\end{aligned}
\end{equation}

For the objective function, let $\Vtctotrg:\bbR^{\numvtc}\to\bbR^{\numtrg}: \Psi\mapsto\Psibar,\quad\Psibar(\trgj) = \vtctotrg(\{\rho(\vtci):\vtci\in\trgj\})$ average the density values on each triangle. 
Below, we list some typical choices of $\vtctotrg:=\vtctotrg\left(\left\{\rho(\vtc_1),\rho(\vtc_2),\rho(\vtc_3)\right\}\right)$.
\begin{itemize}
\item[(i)] Arithmetic mean: $$\omega:=\third\sum_i\rho(\vtci);$$ 
\item[(ii)] Geometric mean: $$\vtctotrg:=\left(\prod_i\rho(\vtci)\right)^{\third};$$
\item[(iii)] Harmonic mean: $$\vtctotrg:=3\left(\sum_i\frac{1}{\rho(\vtci)}\right)^{-1}.$$
\end{itemize}
We remark that these choices of average functions are useful in defining the related discrete mean-field variational problems. They connect with the gradient flow studies of Markov processes on discrete states. See related studies in \cite{MAAS20112250}. For simplicity, we select the arithmetic mean (i) in this work.  

We evaluate the dynamic cost on triangles $\disctcostfuncctl:\discttgsp\to\bbR$, and the interaction and terminal cost on vertices $\disctCostevo:\prob(\disctmfd)\to\bbR$, $\disctCostterm:\prob(\disctmfd)\to\bbR$.
With suitable choice of triangular mesh and discrete cost functions, the continuous cost is approximated by 
\begin{equation}
\begin{aligned}
    \disctCost(\Rho,M):=&
    \int_0^1\sum_{j=1}^{\numtrg} \areatrgj \Rhobar(\trgj,t)\disctcostfuncctl\left(\trgj,\frac{M(\trgj,t)}{\Rhobar(\trgj,t)}\right)\dd t
    +\int_0^1 \disctCostevo(\Rho(\cdot,t))\dd t
    +\disctCostterm(\Rho(\cdot,1)).\\
\end{aligned}
\end{equation}
where $\Rhobar(\cdot,t)=\Vtctotrg(\Rho(\cdot,t))$.
The semi-discrete formulation of \eqref{eq: mfmfg opt} on triangular mesh $\mfd$ is then
\begin{equation}
\begin{aligned}
    \min_{\Rho,M} \quad& \disctCost(\Rho,M)\\ 
    \text{ subject to } &\quad \frac{\dd}{\dd t}\Rho(\vtci,t) + \left(\divg M\right)(\vtci,t) = 0,\Rho(\cdot,0)=\Rho_0,\\
    &\quad  \Rho\in C([0,1];\prob(\disctmfd)), M\in C([0,1];\Gamma(\discttgsp)).
\end{aligned}
\label{eq: mfmfg opt semi-disct}
\end{equation}

\begin{remark}
Recall that in the continuous setting, we show that the local minimizer of the optimization problem \eqref{eq: mfmfg opt} solves the PDE system \eqref{eq: mfmfg kkt} with $\costfuncevo(\vx,\rho(\cdot,t))=\frac{\delta\Costevo(\rho(\cdot,t))}{\delta\rho}(\vx)$ and $\costfuncterm(\vx,\rho(\cdot,1))=\frac{\delta\Costterm(\rho(\cdot,1))}{\delta\rho}(\vx)$.
Since the constraint of \eqref{eq: mfmfg opt semi-disct} remains linear, the local minimizer of this semi-discrete problem also solves a KKT-based ODE system
\begin{equation}\left\{
\begin{aligned}
    &-\frac{\dd}{\dd t}\Phi(\vtci,t) 
    + \sumtrg \frac{\areatrgj}{\areavtci}\frac{\partial\Rhobar(\trgj,t)}{\partial\Rho(\vtci,t)}\discthmt\left( \trgj,(\grad\Phi)(\trgj,t) \right) 
    \leq \frac{1}{\areavtci}\partial_{\Rho(\vtci)}\disctCostevo(\Rho(\cdot,t)),\\
    &\frac{\dd}{\dd t}\Rho(\vtci,t) +\left(\divg M\right)(\vtci,t) = 0,\quad M(\trgj,t) = -\Rhobar(\trgj,t)\partial_{Q}\discthmt(\trgj,(\grad\Phi)(\trgj,t)),\\
    &\Phi(\vtci,1) \leq\frac{1}{\areavtci}\partial_{\Rho(\vtci)}\disctCostterm(\Rho(\cdot,1)),\quad\Rho(\cdot,0)=\Rho_0.
\end{aligned}\right.
\label{eq: mfmfg kkt semi-disct}
\end{equation}
where 
$$\discthmt:\discttgsp\to\bbR, \discthmt(\trgj,Q) = \sup_{U(\trgj)\in \spn{\trgj}}\{-\disctcostfuncctl(\trgj,U(\trgj))-\langle Q,U(\trgj)\rangle\}.$$
Note that $\displaystyle\sumtrg \frac{\areatrgj}{\areavtci}\frac{\partial\Rhobar(\trgj,t)}{\partial\Rho(\vtci,t)}\discthmt\left( \trgj,(\grad\Phi)(\trgj,t) \right)$ is an estimation of $\hmtctl(\vtci,\nabla_{\mfd}\phi(\vtci,t))$. 
If in addition, $\partial_{\Rho(\vtci)}\disctCostevo(\Rho(\cdot,t))/\areavtci$ and $\partial_{\Rho(\vtci)}\disctCostterm(\Rho(\cdot,1))/\areavtci$ approximate the costs $\costfuncevo(\vx)$ and $\costfuncterm(\vx)$ evaluated at $\vx=\vtci$, then \eqref{eq: mfmfg kkt semi-disct} is a semi-discrete formulation of \eqref{eq: mfmfg kkt}.

\end{remark}

Before discretizing the time interval, we discuss the semi-discrete formulations of examples \ref{eg: vmfg local} and \ref{eg: vmfg nonlocal}. 

\begin{example}[Quadratic dynamic cost with local interaction]
\label{eg: vmfg semi-disct local}
This is a discretized counterpart of example \ref{eg: vmfg local}.  
Recall that in the continuous domain, 
Since we choose the surface metric as the induced metric, then for $U(\trgj)\in\spn{\trgj}$, the dynamic cost is provided as
\begin{equation}
    \disctcostfuncctl(\trgj,U(\trgj)) = \half\sum_{d=1}^3(U^d(\trgj))^2=\half\left\|U(\trgj)\right\|_{\bbR^3}^2,
\end{equation}
And we approximate the interaction and terminal costs $\Costevo,\Costterm$ by
\begin{equation}
\begin{aligned}
    &\disctCostevo(\Rho(\cdot,t)) = \sum_{i=1}^{\numvtc}\areavtci\Rho(\vtci,t)\log(\Rho(\vtci,t)),\\
    &\disctCostterm(\Rho(\cdot,1)) =\sum_{i=1}^{\numvtc}\areavtci \Rho(\vtci,1)\log\left(\frac{\Rho(\vtci,1)}{\Rho_1(\vtci)}\right).
\end{aligned}
\end{equation}
With $\vtctotrg$ being the arithmetic average, if the optimizer $\Rho>0$, then the discrete KKT system is
\begin{equation}\left\{
\begin{aligned}
    &-\frac{\dd}{\dd t}\Phi(\vtci,t) 
    + \sumtrg \frac{\areatrgj}{3\areavtci}\half\left\|(\grad\Phi)(\trgj,t) \right\|_{\bbR^3}^2  
    = \log(\Rho(\vtci,t))+1,\\
    &\frac{\dd}{\dd t}\Rho(\vtci,t) + \left(\divg M\right)(\vtci,t) = 0,\quad M(\trgj,t) = -\Rhobar(\trgj,t)(\grad\Phi)(\trgj,t),\\
    &\Phi(\vtci,1) =\log\left(\frac{\Rho(\vtci,1)}{\Rho_1(\vtci)}\right),\quad\Rho(\cdot,0)=\Rho_0.
\end{aligned}\right.
\end{equation}

\end{example}

\begin{example}[Quadratic dynamic cost with non-local interaction cost]
\label{eg: vmfg semi-disct nonlocal}
This is a discretized counterpart of example \ref{eg: vmfg nonlocal}.
On triangular meshes, we use the following dynamic and terminal costs,
\begin{equation}
\begin{aligned}
    &\disctcostfuncctl(\trgj,U(\trgj)) = \half\left\|U(\trgj)\right\|_{\bbR^3}^2,\\
    &\disctCostterm(\Rho(\cdot,1)) = \sum_{i=1}^{\numvtc}\areavtci\half(\Rho(\vtci,1)-\Rho_1(\vtci))^2. \end{aligned}
\end{equation}
And for the discrete interaction cost, we choose $d_{\disctmfd}(\vtc_1,\vtc_2)$ being the length of the shortest path connecting $\vtc_1,\vtc_2$ to approximate the geodesic distance, 
\begin{equation}
d_{\disctmfd}(\vtc_1,\vtc_2)=\min\left\{
\begin{aligned}
    \sum_{k=1}^n\|\vtc_{i_k}-\vtc_{i_{k-1}}\|_{\bbR^3}:&n = 0,1,\cdots,\vtc_{i_0}=\vtc_1,\vtc_{i_n}=\vtc_2,\\
    &\forall k=1,\cdots,n, \text{ there exists }\trg_{j_k} \text{ such that }\vtc_{i_{k-1}},\vtc_{i_k}\in\trg_{j_k}, 
\end{aligned}
\right\}.
\end{equation}
And we define the kernel as $\disctkernel(\vtc_1,\vtc_2)=\mu\exp\left(-\frac{1}{\sigma}d_{\disctmfd}^2(\vtc_1,\vtc_2)\right)$ and the interaction cost as
\begin{equation}
    \disctCostevo(\Rho(\cdot,t)) =\half \sum_{i}^{\numvtc}\sum_{i'}^{\numvtc}\areavtci\Rho(\vtci) \disctkernel(\vtc_i,\vtc_{i'})\area_{\vtc_{i'}}\Rho(\vtc_{i'})=\half\Rho(\cdot,t)^\top\areavtc \disctkernel\areavtc\Rho(\cdot,t).
\end{equation}
Similarly, if $\vtctotrg$ is the arithmetic average and the optimizer $\Rho>0$, the discrete KKT system is
\begin{equation}\left\{
\begin{aligned}
    &-\frac{\dd}{\dd t}\Phi(\vtci,t) 
    + \sumtrg \frac{\areatrgj}{3\areavtci}\half\left\|(\grad\Phi)(\trgj,t) \right\|_{\bbR^3}^2  
    =  \sum_{i'}\disctkernel(\vtc_i,\vtc_{i'})\area_{\vtc_{i'}}\Rho(\vtc_{i'}),\\
    &\frac{\dd}{\dd t}\Rho(\vtci,t) + \left(\divg M\right)(\vtci,t) = 0,\quad M(\trgj,t) = -\Rhobar(\trgj,t)(\grad\Phi)(\trgj,t),\\
    &\Phi(\vtci,1) =\Rho(\vtci,1)-\Rho_1(\vtci),\quad\Rho(\cdot,0)=\Rho_0.
\end{aligned}\right.
\end{equation}
\end{example}

\subsection{Time discretization}
\label{subsec: disct time}

To numerically solve \eqref{eq: mfmfg opt semi-disct}, we fully discretize the problem by dividing the time interval $[0,1]$ into $n$ segments and let $\tk=\frac{k}{n}$.
Now we consider the density on central time steps $\Rho=\{\Rho(\cdot,\tk)\}_{k=1,\cdots,n}\in(\prob(\disctmfd))^n$ and the flux on staggered time steps $M=\{M(\cdot,\tkm)\}_{k=1,\cdots,n}\in(\Gamma(\discttgsp))^n$.

Let the time differential operator be
\begin{equation}
    (\Dt\Rho)(\cdot,\tkm):=\begin{cases}
    \frac{1}{1/n}(\Rho(\cdot,\tk)-\Rho(\cdot,\tkmm)),& \quad k=2,\cdots,n,\\
    \frac{1}{1/n}(\Rho(\cdot,t_1)-\Rho_0(\cdot)),&\quad k=1.
    \end{cases}
\end{equation}
Then the discrete constraint set $\disctcstr(\Rho_0)$ is
\begin{equation}
\disctcstr(\Rho_0):=\left\{
\begin{aligned}
    (\Rho,M):  &~(\Dt\Rho)(\vtci,\tkm) + (\divg M)(\vtci,\tkm) = 0,\forall\vtci\in\vtc,k=1,\cdots,n\\
    & \quad\Rho\in\prob(\disctmfd))^n, M\in (\Gamma(\discttgsp))^n 
\end{aligned}\right\}.
\end{equation}
Additionally, letting
\begin{equation}
\begin{aligned}
    \Rhobar(\cdot,\tkm)&:=\half\Vtctotrg(\Rho(\cdot,\tk)) + \half\Vtctotrg(\Rho(\cdot,\tkmm)),\quad k=1,\cdots,n\\
    \disctCost(\Rho,M)&:=\frac{1}{n}\sum_{k=1}^{n}\sum_{j=1}^{\numtrg} \areatrgj \Rhobar(\trgj,\tkm)\costfuncctl\left(\trgj,\frac{M(\trgj,\tkm)}{\Rhobar(\trgj,\tkm)}\right) 
    + \frac{1}{n}\sum_{k=1}^{n-1}\disctCostevo(\Rho(\cdot,\tk))
    + \disctCostterm(\Rho(\cdot,\tn)).
\end{aligned}
\end{equation}
we formulate the discrete optimization problem as 
\begin{equation}
\begin{aligned}
    &\min_{\Rho,M}
    \disctCost(\Rho,M) + \chi_{\disctcstr(\rho_0)}(\Rho,M).\\
\end{aligned}
\label{eq: mfmfg opt disct}
\end{equation}
Here $\chi$ is the indicator function
$\chi_\cstr(\vx)=\left\{\begin{array}{cc}
    0, & \vx \in\cstr\\
    +\infty, & \vx\not\in\cstr
    \end{array}\right.$ of a convex set $\cstr$. 

In the next section, we focus on solving the optimization problem \eqref{eq: mfmfg opt disct}.

\section{Algorithm for Solving Variational MFGs on Triangular Meshes}
\label{sec: alg}

In this section, we adapt the fast algorithm proposed in \cite{yu2021fast} to solve the discretized potential mean-field game \eqref{eq: mfmfg opt disct}. This algorithm is based on a proximal gradient method (PGD)~\cite{rockafellar1970convex,beck2009fast}.

To solve \eqref{eq: mfmfg opt disct}, we conduct gradient descent on the smooth component $\disctCost$ of the objective function and proximal descent on the non-smooth component $\chi_{\disctcstr(\Rho_0)}$.
The gradient descent step is trivially
\begin{equation}
    \left(\Rho^\nitp,M^\nitp\right) 
    = \left(\Rho^\nit,M^\nit\right) 
    - \eta^\nit\nabla_{\Rho,M}\disctCost \left(\Rho^\nit,M^\nit\right) 
\end{equation}
with stepsize $\eta^\nit$.
The proximal descent is exactly the projection to $\disctcstr(\Rho_0)$. Let the inner product in discrete spaces be 
\begin{equation*}
\begin{aligned}
    &\Psi_1,\Psi_2\in(\prob(\disctmfd))^n,
    &\langle \Psi_1,\Psi_2\rangle_{\vtc,t} := \frac{1}{n}\sum_{k=1}^n \langle \Psi_1(\cdot,\tk),\Psi_2(\cdot,\tk) \rangle_{\vtc},\\
    &U_1,U_2\in(\Gamma(\discttgsp))^n,
    &\langle U_1,U_2\rangle_{\trg,t} := \frac{1}{n}\sum_{k=1}^{n}\langle U_1(\cdot,\tkm),U_2(\cdot,\tkm) \rangle_{\trg},
\end{aligned}
\end{equation*}
Then 
\begin{equation}
    (\Rho^{\nitpp},M^{\nitpp})=\proj_{\disctcstr(\Rho_0)}(\Rho^{\nitp},M^{\nitp}):=\argmin_{(\Rho,M)\in\disctcstr(\Rho_0)} 
    \half\left\|\Rho-\Rho^{\nitp}\right\|^2_{\vtc,t} +
    \half\left\|M-M^{\nitp}\right\|^2_{\trg,t}
\label{eq: mfmfg disct proj}
\end{equation}
To solve the optimization problem \eqref{eq: mfmfg disct proj}, we introduce a dual variable $\Psi=\{\Psi(\cdot,\tkm)\}_{k=1,\cdots,n}\in(\prob(\disctmfd))^n$ on vertices and the staggered time steps. The Lagrangian is therefore
\begin{equation}
\begin{aligned}
    \Lag(\Rho,M,\Psi):=&\half\left\|\Rho-\Rho^{\nitp}\right\|^2_{\vtc,t} 
    +\half\left\|M-M^{\nitp}\right\|^2_{\trg,t}  \\
    &+\frac{1}{n}\sum_{k=1}^{n}\langle\Psi(\cdot,\tkm), (\Dt\Rho)(\cdot,\tkm)+\divg M(\cdot,\tkm) \rangle_{\vtc}.
\end{aligned}
\end{equation}
If we define $\Dt^*$ as
\begin{equation}
    (\Dt^*\Psi)(\cdot,\tk):=
    \begin{cases}
    \frac{1}{1/n}(\Psi(\cdot,\tkm)-\Psi(\cdot,\tkp)), &k=1,2,\cdots,n-1,\\
    \frac{1}{1/n}\Psi(\cdot,t_{n-\half}),& k=n,
    \end{cases}
\end{equation}
then the Lagrangian is also 
\begin{equation}
\begin{aligned}
    \Lag(\Rho,M,\Psi):=&\half\left\|\Rho-\Rho^{\nitp}\right\|^2_{\vtc,t} 
    +\half\left\|M-M^{\nitp}\right\|^2_{\trg,t}  \\
    &+\frac{1}{n}\sum_{k=1}^{n}\langle (\Dt^*\Psi)(\cdot,\tk),\Rho(\cdot,\tk) \rangle_{\vtc}-\langle \Psi(\cdot,t_{\half}),\Rho_0 \rangle_{\vtc}\\
    &+\frac{1}{n}\sum_{k=1}^{n}\langle -\grad\Psi(\cdot,\tkm),M(\cdot,\tkm) \rangle_{\trg}.
\end{aligned}
\end{equation}
Thus the saddle point $(\Rho,M,\Psi)$ satisfies the linear system
\begin{equation}
    (\Dt\Rho)(\cdot,\tkm)+\divg M(\cdot,\tkm)=\Zero,\quad k=1,\cdots,n.
\end{equation}
and
\begin{equation}
\left\{\begin{aligned}
    \Rho(\cdot,\tk) &= \Rho^{\nitp}(\cdot,\tk) - (\Dt^*\Psi)(\cdot,\tk),\quad k=1,\cdots,n\\
    M(\cdot,\tkm) &= M^{\nitp}(\cdot,\tkm) + \grad\Psi(\cdot,\tkm),\quad k=1,\cdots,n,
\end{aligned}\right.
\end{equation}
Note that $\Dt$ is a full rank operator, for $k=1,\cdots,n$, $\Psi(\cdot,\tkm)$ is the unique solution to
\begin{equation}
    (\Dt\Dt^*\Psi)(\cdot,\tkm) - \divg\grad \Psi(\cdot,\tkm) = 
    (\Dt\Rho^{\nitp})(\cdot,\tkm) + \divg M^{\nitp}(\cdot,\tkm).
\label{eq: mfmfg phi equ}
\end{equation}
Since this linear solver is invariant to the data and the iteration number, in practice, we precompute it to save cost in the main iteration.

We summarize our algorithm in Alg. \ref{alg: pgd for mfmfg}.

\begin{algorithm}
\caption{PGD for MFG on discrete mesh}
\begin{algorithmic}
    \STATE{Parameters}
    $\Rho_0$
    \STATE{Initialization} 
    $\Rho^{(0)}\in(\prob(\disctmfd))^n$, and $M^{(0)}\in(\Gamma(\discttgsp))^n$.
    
    \FOR{$l=0,1,2,\ldots$}
    \STATE{\textbf{gradient descent}} 
    \begin{equation*}
    \left(\Rho^\nitp,M^\nitp\right) 
    = \left(\Rho^\nit,M^\nit\right) 
    - \eta^\nit\nabla_{\Rho,M}\disctCost \left(\Rho^\nit,M^\nit\right) 
    \end{equation*}
    
    \STATE{\textbf{proximal descent}}\quad for $k=1,\cdots,n$, solve $\Psi$ for 
    \begin{equation*}
        (\Dt\Dt^*\Psi)(\cdot,\tkm) - \divg\grad \Psi(\cdot,\tkm) = 
        (\Dt\Rho^{\nitp})(\cdot,\tkm) + \divg M^{\nitp}(\cdot,\tkm).
    \end{equation*}
    and conduct
    \begin{equation*}
    \left\{\begin{aligned}
        \Rho^\nitpp(\cdot,\tk) &= \Rho^{\nitp}(\cdot,\tk) - (\Dt^*\Psi)(\cdot,\tk),\quad k=1,\cdots,n\\
        M^\nitpp(\cdot,\tkm) &= M^{\nitp}(\cdot,\tkm) + \grad\Psi(\cdot,\tkm),\quad k=1,\cdots,n,
    \end{aligned}\right.
    \end{equation*}
    
    \ENDFOR
\end{algorithmic}
\label{alg: pgd for mfmfg}
\end{algorithm}

\section{Numerical Examples}
\label{sec: num res}

In this section, we conduct various experiments to show the effectiveness and flexibility of our mean-field game models on manifolds and the proposed numerical method.
We provide numerical results on different manifolds. Most of these manifolds are non-Euclidean, thus conventional settings of mean-field games cannot handle them.
In all of our experiments, we choose the induced metric for the manifold geometry, the quadratic dynamic cost $\costfuncctl\left(\trgj,U(\trgj,t)\right):=\half\sum_d \left(U^d(\trgj,t)\right)^2$,
and the arithmetic average $\vtctotrg$.
The interaction and terminal cost terms vary from examples and will be specified later.
All of our numerical experiments are implemented in Matlab on a PC with an Intel(R) i7-8550U 1.80GHz CPU  and 16 GB memory. 

\subsection{MFGs with local interactions}

In this part, both $\Costevo$ and $\Costterm$ take local forms for all experiments.

\begin{figure}[ht]
\flushright
\includegraphics[scale=1]{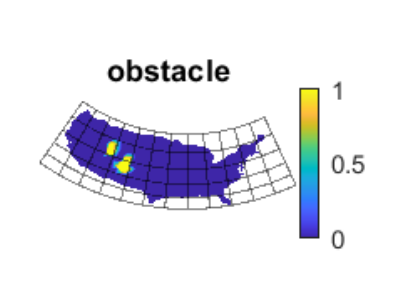}
\includegraphics[scale=1]{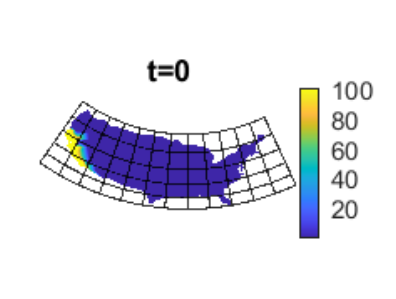}
\includegraphics[scale=1]{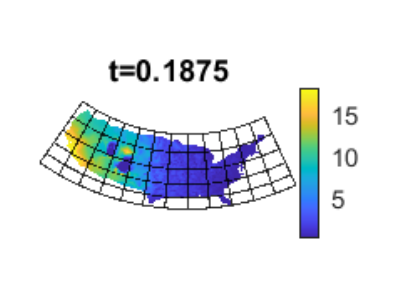}
\includegraphics[scale=1]{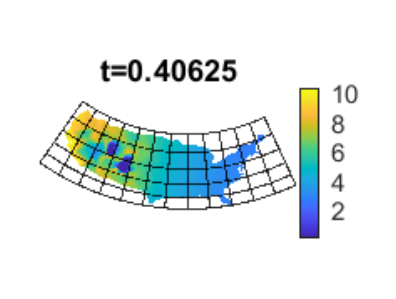}
\includegraphics[scale=1]{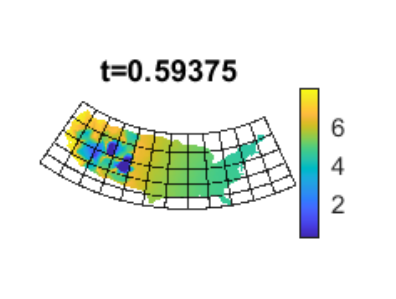}
\includegraphics[scale=1]{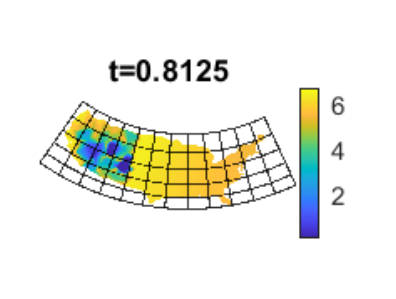}
\includegraphics[scale=1]{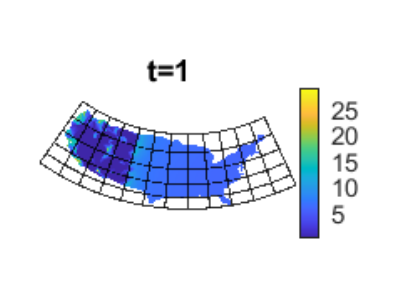}
\caption{Illustration of obstacle indicator $B$ (column 1) and snapshots of a MFG on the US map (column 2-4). 
}
\label{fig: usmap result}
\end{figure}

\paragraph{The U.S. map based triangular mesh}
We first consider a U.S. map based triangular mesh, which is the discretization of a subdomain on a spherical manifold.
Assume that there are two obstacles on the map and it takes extra efforts for masses (agents) to pass through the obstacle region $\disctcstr_{B}\subset\disctmfd$.
We define $B:\disctmfd\to\bbR$ be the piece-wise linear indicator of the obstacle with
$B(\vtci)=\begin{cases}
1,&\vtci\in\disctcstr_B,\\
0,&\vtci\in\disctcstr_B
\end{cases}$ (see Figure \ref{fig: usmap result} obstacle). 

We pick the initial density $\Rho_0$ showed in Figure \ref{fig: usmap result} $t=0$. The mass concentrates in California.
We let the mass to move freely during the time interval.
To reflect the impact of the obstacle, we choose the interaction cost $\disctCostevo(\Rho(\cdot,t)) =50\sum_{i=1}^{\numvtc}\areavtci \Rho(\vtci,t)B(\vtci)$.
We also encourage the mass to stop in the central and eastern part at the end.
To achieve this, we define the terminal cost $\disctCostterm(\Rho(\cdot,t)) = \frac{1}{10}\sum_{i=1}^{\numvtc}\areavtci\Rho(\vtci,t_n)B_{T}(\vtci)$ with
$B_T(\vtci)=\begin{cases}
1,& 105^\circ \text{W} \leq\text{the longitude of } \vtci  \leq 130^\circ\text{W},\\
0,&\text{otherwise}.
\end{cases}$

Our numerical results in Figure \ref{fig: usmap result} show that the density in the obstacle region remains low. This means the mass circumvent those areas very well.
In addition, at $t=1$, the density in the central and eastern areas is generally denser than that in the western, which meets our expectation.

\paragraph{``8''-shape with obstacles}

\begin{figure}[h]
\centering
\subfigure[Illustration of the manifold, obstacle, initial density and terminal density]{
\label{fig: eightobs illus}
\includegraphics[scale=1]{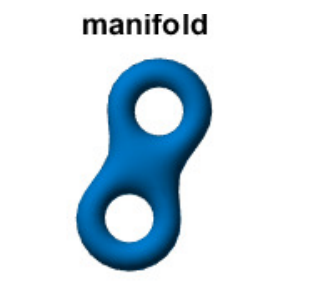}
\includegraphics[scale=1]{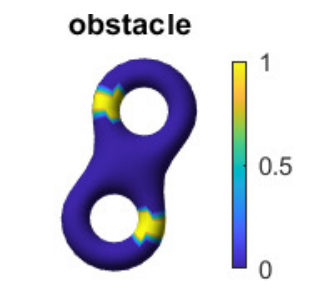}
\includegraphics[scale=1]{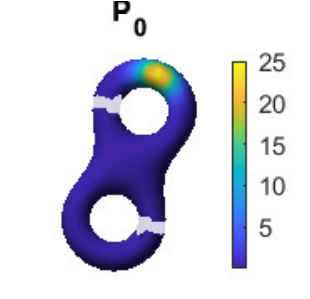}
\includegraphics[scale=1]{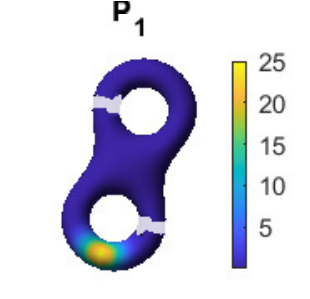}
}
\centering
\subfigure[Snapshots of the evolution]{
\label{fig: eightobs shots}
\begin{minipage}[b]{\textwidth}
\includegraphics[scale=1]{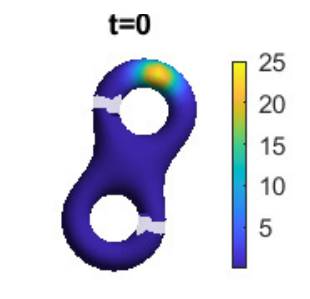}
\includegraphics[scale=1]{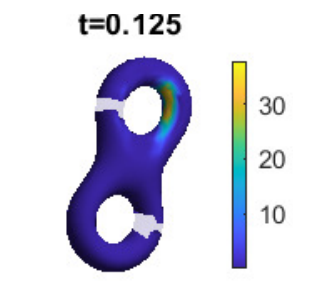}
\includegraphics[scale=1]{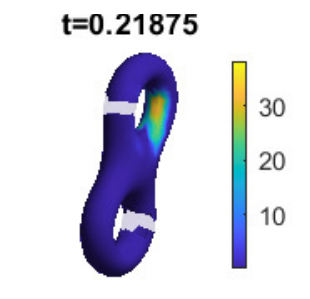}
\includegraphics[scale=1]{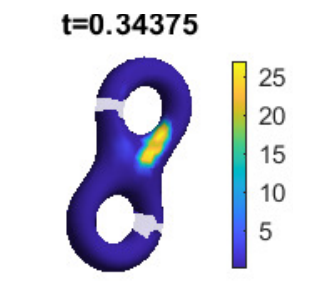}
\includegraphics[scale=1]{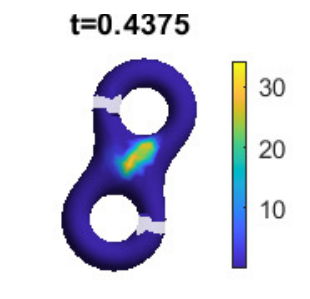}
\includegraphics[scale=1]{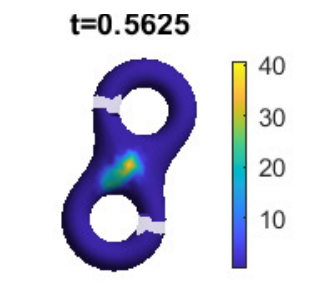}
\includegraphics[scale=1]{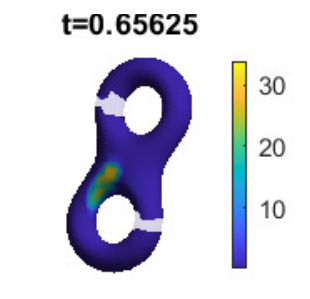}
\includegraphics[scale=1]{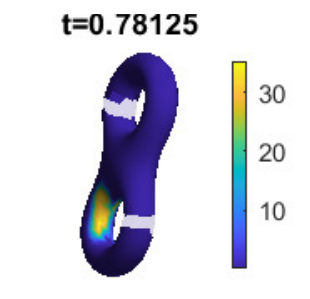}
\includegraphics[scale=1]{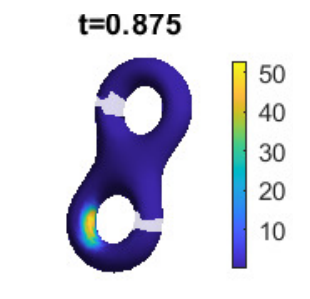}
\includegraphics[scale=1]{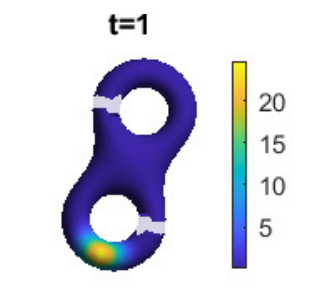}
\end{minipage}
}
\caption{Illustration of the model and snapshots of the density evolution.}
\label{fig: eightobs}
\end{figure}

In this example, we demonstrate that our model and algorithm can successfully handle manifolds with complicated topology. 
We consider a  ``8'' shape surface. 
Similarly as before, we assume there are obstacles on the manifold and the indicator of the obstacle $B:\disctmfd\to\bbR$ is shown in Figure \ref{fig: eightobs illus}.
In plots afterwards, we indicate the obstacle region with a different transparency.
We pick the initial density $\Rho_0$ aggregating on the one end of ``8'' and the desired terminal density $\Rho_1$ on the other end (Figure \ref{fig: eightobs illus}).
For the interaction cost, we still choose $\disctCostevo(\Rho(\cdot,t)) = 50\sum_{i=1}^{\numvtc}\areavtci \Rho(\vtci,t)B(\vtci)$ to avoid obstacle.
And we write the terminal cost as $\disctCostterm(\Rho(\cdot,t)) = 5\sum_{i=1}^{\numvtc}\areavtci\left(\Rho(\vtci,1)-\Rho_1(\vtci)\right)^2$ to push the terminal density $\Rho(\cdot,1)$ to the desired $\Rho_1$.

We list the snapshots of resulting density evolution in Figure \ref{fig: eightobs shots}. 
These results show that the mass produced by our model successfully circumvents the obstacle on this genus-2 manifold. Additionally, the terminal density mainly aggregated in the support of $\Rho_1$ as we expect.

\paragraph{Irregular Euclidean domain}

Besides introducing $B$ to impose a soft constraint of the obstacle, the general setup on manifolds enables us to have a different implementation for the hard obstacle constraints.

\begin{figure}[h]
\centering
\subfigure[MFG with vanilla $\disctCostevo_v(\Rho(\cdot,t)) = 0$]{
\label{fig: eucobs vanila}
\includegraphics[scale=1]{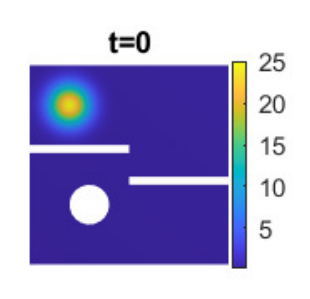}
\includegraphics[scale=1]{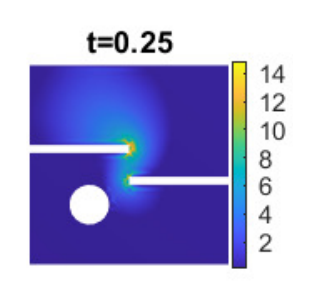}
\includegraphics[scale=1]{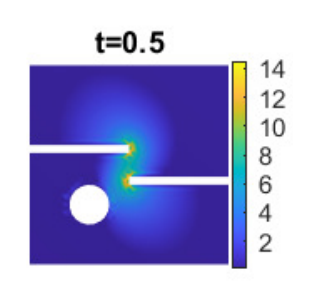}
\includegraphics[scale=1]{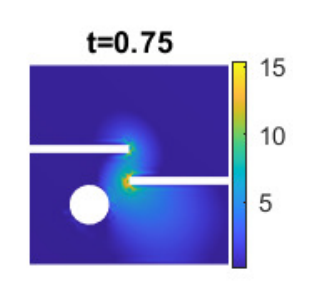}
\includegraphics[scale=1]{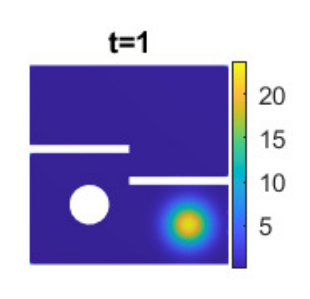}
}
\subfigure[MFG with disperse $\disctCostevo_d(\Rho(\cdot,t)) = \sum_{i=1}^{\numvtc}\areavtci\Rho(\vtci,t)\log(\Rho(\vtci,t))$]{
\label{fig: eucobs rhologrho}
\includegraphics[scale=1]{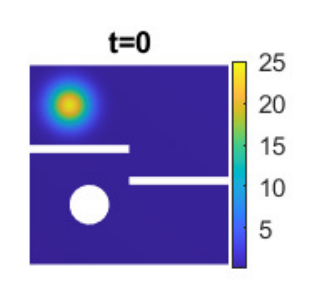}
\includegraphics[scale=1]{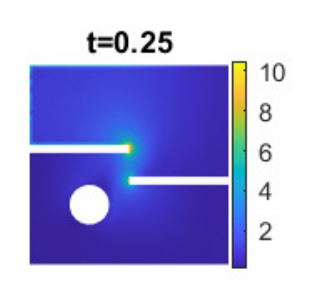}
\includegraphics[scale=1]{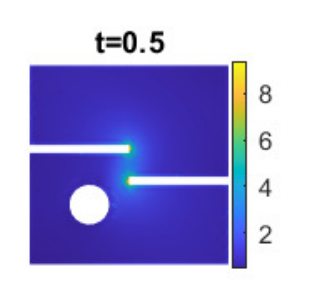}
\includegraphics[scale=1]{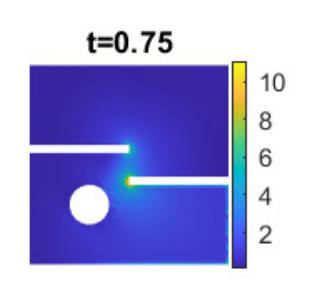}
\includegraphics[scale=1]{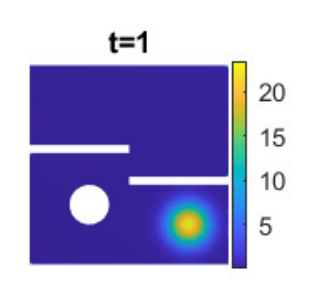}
}
\caption{Snapshots of MFGs with different interactions on constrained Euclidean space.}
\label{fig: eucobs compare}
\end{figure}

For example, consider an irregular Euclidean domain showed in Figure \ref{fig: eucobs compare} with white regions punctured. Instead of handling the complicate boundary conditions when conducting mean-field game problems using conventional methods in Euclidean spaces, 
We view the region as a two-dimensional manifold and use a triangular mesh to approximate it. Then we directly apply our algorithm to the mesh without concerning the shape complexity.

The initial density $\Rho_0$ and desired terminal density $\Rho_1$ are approximations of two Gaussian distributions. 
And we choose terminal cost as $\disctCostterm(\Rho(\cdot,t)) = 10\sum_{i=1}^{\numvtc}\areavtci\Rho(\vtci,1)\log\left(\frac{\Rho(\vtci,1)}{\Rho_1(\vtci)}\right)$ to push $\Rho(\cdot,1)$ to $\Rho_1$.
In Figure \ref{fig: eucobs compare}, we compare the results with a vanilla interaction cost $\disctCostevo_v(\Rho(\cdot,t)) = 0$ (Figure \ref{fig: eucobs vanila}) and with a disperse cost $\disctCostevo_d(\Rho(\cdot,t)) = \sum_{i=1}^{\numvtc}\areavtci\Rho(\vtci,t)\log(\Rho(\vtci,t))$ (Figure \ref{fig: eucobs rhologrho}).
We see that with $\disctCostevo_d$, the mass is prone to segregate during the evolution.
To understand this, we refer to the original game description. 
With $\disctCostevo_d$, we actually solve a mean-field game with $\costfuncevo(\vx,\rho(\cdot,t))=\log(\rho(\vx,t))+1$.
To reduce the cost $\Jfunc$, agents prefer locations with lower $\costfuncevo(\vx,\rho(\cdot,t))$, i.e. lower density $\rho(\vx)$.

\paragraph{Homer surface}

As the last example with local cost, we work with the surface of homer.
We pick the initial density $\Rho_0$ concentrating on the belly and the desired terminal density $\Rho_1$ on the end of hands and feet.
Fixing the terminal cost $\disctCostterm(\Rho(\cdot,1))=\frac{5}{2}\sum_{i=1}^{\numvtc}\areavtci(\Rho(\vtci,1)-\Rho_1(\vtci))^2$, we compare the vanilla interaction cost $\disctCostevo_v(\Rho(\cdot,t))=0$ and congested $\disctCostevo_c(\Rho(\cdot,t))=\frac{1}{10}\sum_{i=1}^{\numvtc}\areavtci\sqrt{\Rho(\vtci,t)+10^{-4}}$.
The choice of $\disctCostevo_c$ actually corresponds to the mean-field game with $\costfuncevo(\rho(\cdot,t))=\frac{1}{10\sqrt{\rho(\vx,t)+10^4}}$.
To reduce the cost during evolution, the agents tend to aggregate for a larger density value.

\begin{figure}[H]
\centering
\subfigure[MFG with vanilla $\disctCostevo_v(\Rho(\cdot,t)) = 0$]{
\label{fig: homer vanila}
\includegraphics[scale=1]{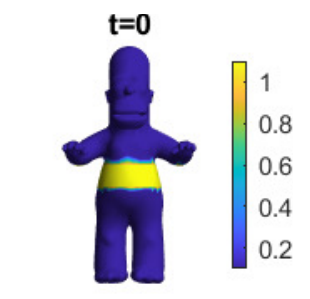}
\includegraphics[scale=1]{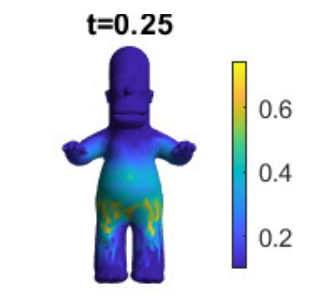}
\includegraphics[scale=1]{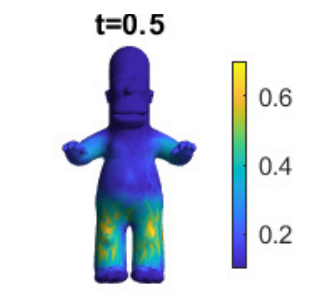}
\includegraphics[scale=1]{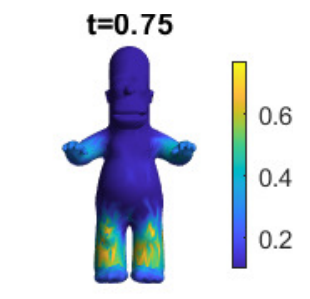}
\includegraphics[scale=1]{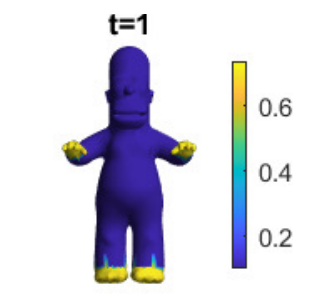}
}
\subfigure[MFG with congested $\disctCostevo_c(\Rho(\cdot,t))=\frac{1}{10}\sum_{i=1}^{\numvtc}\areavtci\sqrt{\Rho(\vtci,t)+10^{-4}}$]{
\label{fig: homer sqrtrho}
\includegraphics[scale=1]{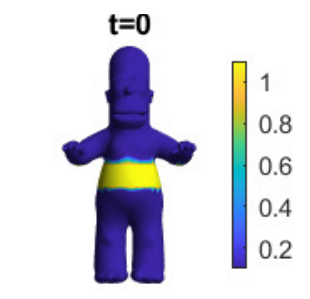}
\includegraphics[scale=1]{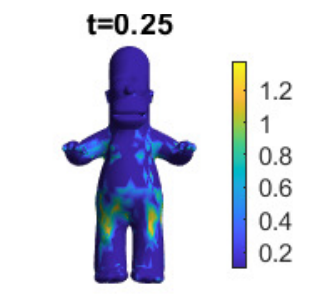}
\includegraphics[scale=1]{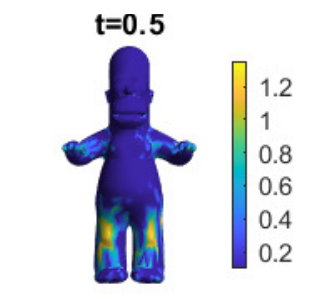}
\includegraphics[scale=1]{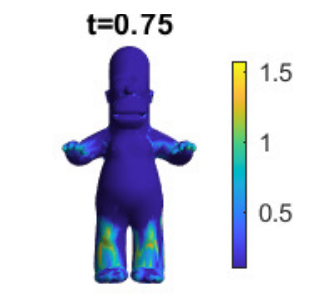}
\includegraphics[scale=1]{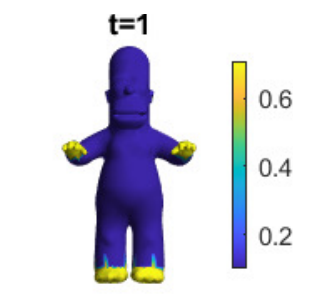}
}
\caption{Snapshots of a MFG on the homer surface.}
\label{fig: homer result}
\end{figure}

We solve the games with $\disctCostevo_v,\disctCostevo_c$ to obtain the local minimizers $(\Rho_v,M_v),(\Rho_c,M_c)$ and report the corresponding dynamic cost, terminal cost and value of $\disctCostevo_c$ in Table \ref{tab: homer cost compare}.
We also show and compare $\Rho_v(\cdot,t),\Rho_c(\cdot,t)$ at several time steps in Figure \ref{fig: homer result}.
The costs in Table \ref{tab: homer cost compare} show that our algorithm effectively reduce the interaction cost.
From the snapshots, we observe that with the congested interaction cost, the mass move in a more compact manner.

\begin{table}[thb]
\centering
\caption{Comparison of dynamic, interaction and terminal costs for experiments on the homer surface.}
\begin{tabular}{l|ccc}
\toprule
          & dynamic cost & $\frac{1}{n}\sum_{k=1}^{n-1}\disctCostevo_c(\Rho(\cdot,\tk))$  & terminal cost \\ \midrule
vanilla   & 0.0079  & 0.0393    & $5.5\times10^{-5}$   \\
congested & 0.0084  & 0.0378    & $8.2\times10^{-5}$   \\
\bottomrule
\end{tabular}
\label{tab: homer cost compare}
\end{table}

\subsection{MFGs with non-local interactions}

In this part, we show some mean-field games with non-local interaction cost. 

\paragraph{The unit Sphere}

\begin{figure}[h]
\centering
\subfigure[MFG with vanilla $\disctCostevo_v(\Rho(\cdot,t)) = 0$.]{
\label{fig: sphere vanila}
\includegraphics[scale=1]{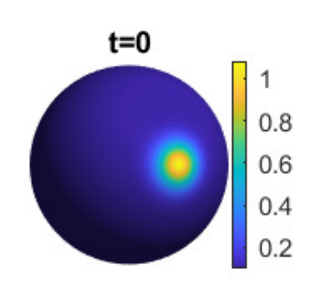}
\includegraphics[scale=1]{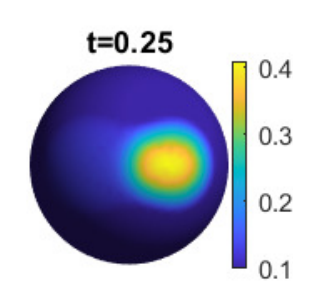}
\includegraphics[scale=1]{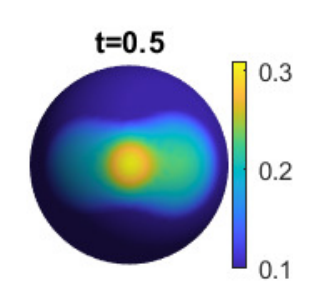}
\includegraphics[scale=1]{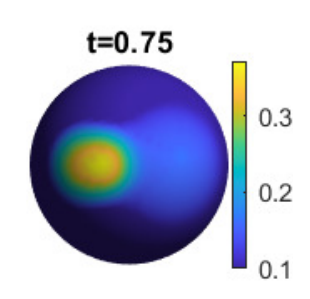}
\includegraphics[scale=1]{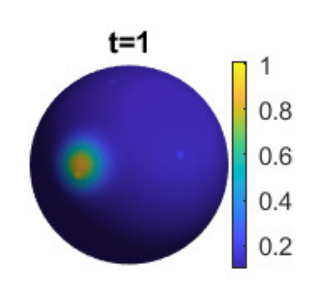}
}
\subfigure[MFG with non-local $\disctCostevo_n(\Rho(\cdot,t)) =25 \sum_{i=1}^{\numvtc}\sum_{i'}\areavtci\Rho(\vtci)\disctkernel(\vtc_i,\vtc_{i'})\area_{\vtc_{i'}}\Rho(\vtc_{i'})$.]{
\label{fig: sphere nonlocal}
\includegraphics[scale=1]{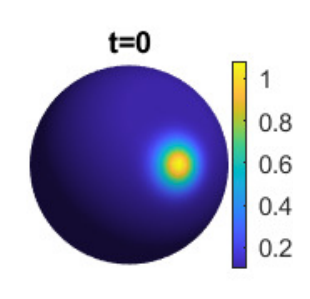}
\includegraphics[scale=1]{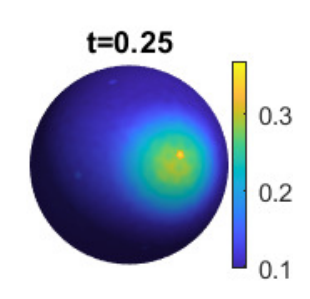}
\includegraphics[scale=1]{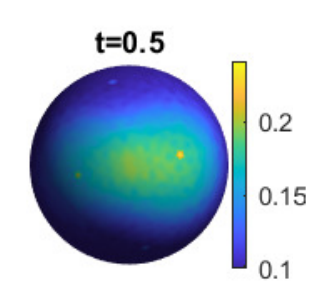}
\includegraphics[scale=1]{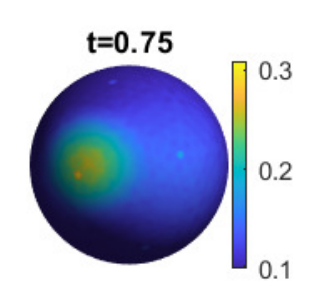}
\includegraphics[scale=1]{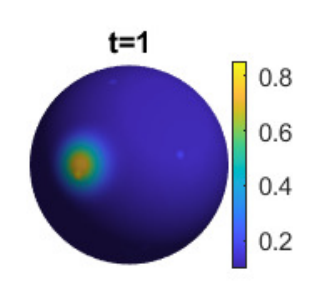}
}
\caption{Snapshots of MFGs with different interactions on sphere.}
\label{fig: sphere compare}
\end{figure}

In this example, we work on the triangular mesh of the unit sphere in three-dimensional space.
The initial density $\Rho_0$ and desired terminal density $\Rho_1$ are spherical Gaussian.
Again, we use the terminal cost
$\disctCostterm(\Rho(\cdot,t)) = 0.5\sum_{i=1}^{\numvtc}\areavtci\Rho(\vtci,1)\log\left(\frac{\Rho(\vtci,1)}{\Rho_1(\vtci)}\right).$ 
We then compute the game with vanilla interaction cost $\disctCostevo_v(\Rho(\cdot,t)) =0$ and non-local $\disctCostevo_n(\Rho(\cdot,t)) =25 \sum_{i=1}^{\numvtc}\sum_{i'}\areavtci\Rho(\vtci)\disctkernel(\vtc_i,\vtc_{i'})\area_{\vtc_{i'}}\Rho(\vtc_{i'})$.
The kernel is defined as
$$\disctkernel(\vtci,\vtc_{i'}) = \exp\left(-(\arccos{\vtci^\top\vtc_{i'}})^2/\sigma^2\right).$$
Here $\sigma=0.1$ and $\vtci^\top\vtc_{i'}$ is the inner product of the two vectors in Euclidean space and $\arccos{\vtci^\top\vtc_{i'}}$ is the geodesic distance between $\vtci$ and $\vtc_{i'}$ on sphere. 
We use the ground truth geodesic distance for simplicity. One can also compute the shortest path on the mesh and store it when pre-processing the manifold.
We conduct the quantitative and snapshot comparison in Table \ref{tab: sphere cost compare} and in Figure \ref{fig: sphere compare}.
The table shows our algorithm effectively leverage the dynamic, interaction and terminal cost when taking the non-local cost $\disctCostevo_n$.
And the comparison in Figure \ref{fig: sphere compare} clearly illustrates that the non-local cost $\disctCostevo_n$ encourages the dispersion of mass.

\begin{table}[thb]
\centering
\caption{Comparison of dynamic, interaction and terminal costs for experiments on sphere.}
\begin{tabular}{l|ccc}
\toprule
          & dynamic cost & $\frac{1}{n}\sum_{k=1}^{n-1}\disctCostevo_n(\Rho(\cdot,\tk))$ & terminal cost \\ \midrule
vanilla   & 0.0267  & 0.1220    & 0.0023   \\
non-local & 0.0292  & 0.1148    & 0.0039   \\
\bottomrule
\end{tabular}
\label{tab: sphere cost compare}
\end{table}

\paragraph{Kitten}

\begin{figure}[h]
\centering
\subfigure[(vanilla) MFG with $\disctCostevo_v(\Rho(\cdot,t)) = 0$]{
\label{fig: kitten vanila}
\includegraphics[scale=1]{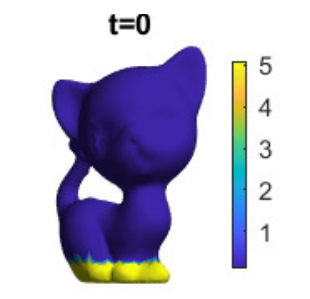}
\includegraphics[scale=1]{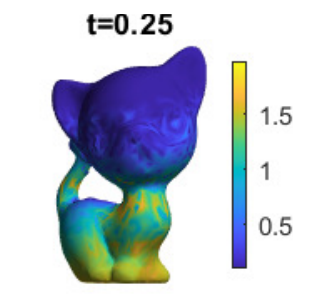}
\includegraphics[scale=1]{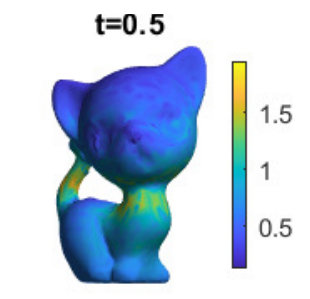}
\includegraphics[scale=1]{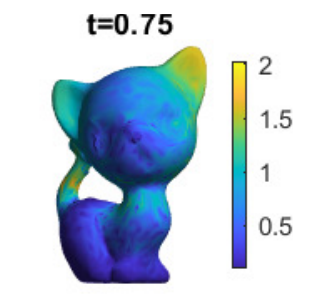}
\includegraphics[scale=1]{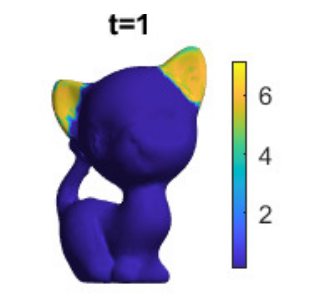}
}
\subfigure[(non-local) MFG with $\disctCostevo_n(\Rho(\cdot,t)) =\half \sum_{i=1}^{\numvtc}\sum_{i'}\areavtci\Rho(\vtci,t)\disctkernel(\vtc_i,\vtc_{i'})\area_{\vtc_{i'}}\Rho(\vtc_{i'},t)$]{
\label{fig: kitten nonlocal}
\includegraphics[scale=1]{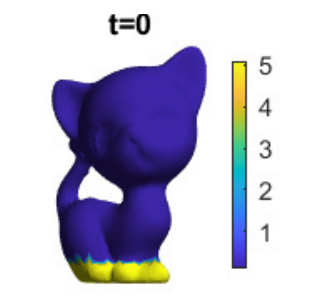}
\includegraphics[scale=1]{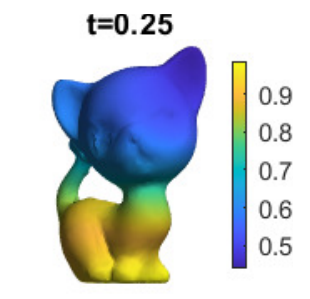}
\includegraphics[scale=1]{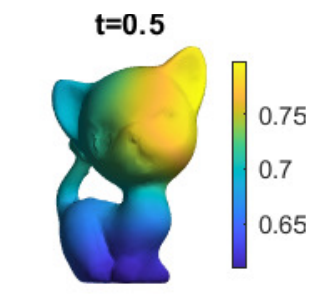}
\includegraphics[scale=1]{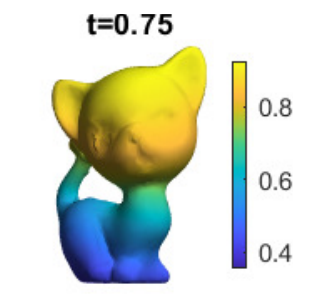}
\includegraphics[scale=1]{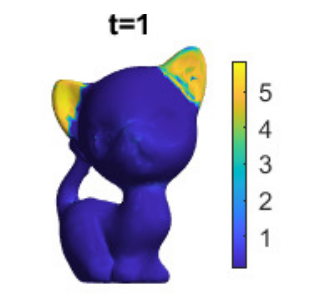}
}
\caption{Snapshots of MFGs with different interactions on ``kitten''.}
\label{fig: kitten compare}
\end{figure}

In the last example, we work with the kitten surface (Figure \ref{fig: kitten compare}).
Let the initial density $\Rho_0$ to concentrate on the paws and the desired terminal density on ears.
We take the terminal cost $\disctCostterm(\Rho(\cdot,1)) = \sum_{i=1}^{\numvtc}\areavtci\left(\Rho(\vtci,1)-\Rho_1(\vtci)\right)^2$ to push the mass moving from bottom to top. 
We also compare the non-local interaction cost $\disctCostevo_n(\Rho(\cdot,t)) =\half \sum_{i=1}^{\numvtc}\sum_{i'}\areavtci\Rho(\vtci,t)\disctkernel(\vtc_i,\vtc_{i'})\area_{\vtc_{i'}}\Rho(\vtc_{i'},t)$ with the vanilla $\disctCostevo_v(\Rho(\cdot,t))=0$ in Figure \ref{fig: kitten compare}.
The kernel is chosen as a weighted Laplacian matrix on the triangular mesh
$$\disctkernel(\vtci,\vtc_{i'})=\frac{1}{\areavtci}\frac{1}{\area_{\vtc_{i'}}}
\sum_{d=1}^3\sum_{j=1}^{\numtrg}\areatrgj\disctmtc^d(\trgj,\vtci)\disctmtc^d(\trgj,\vtc_{i'}).$$
With this choice, the interaction cost is exactly 
$$\disctCostevo_n(\Rho(\cdot,t))=\half\sum_{j=1}^{\numtrg}\areatrgj\left\|(\grad\Rho)(\trgj,\tk)\right\|_2^2,$$
and approximates $\half\int_{\mfd}\left\|\nabla_{\mfd}\rho(\vx,t)\right\|_{\mtc(\vx)}^2\ddx$.
To reduce this cost, the density at each time step $\Rho(\cdot,\tk)$ tends to be smooth on the space domain.
The quantitative result in Table \ref{tab: kitten cost compare} shows the value $\frac{1}{n}\sum_{k=1}^{n-1}\disctCostevo_n(\Rho(\cdot,\tk))$ is reduced by adding $\disctCostevo_n$ to the objective function.
And the comparisons of densities and colorbars in Figure \ref{fig: kitten compare} show that with $\disctCostevo_n$ in objective function, at each time step, the density distributes more uniformly on the manifold.

\begin{table}[thb]
\centering
\caption{Comparison of dynamic, interaction and terminal costs for experiments on ``kitten''.}
\begin{tabular}{l|ccc}
\toprule
          & dynamic cost & $\frac{1}{n}\sum_{k=1}^{n-1}\disctCostevo_n(\Rho(\cdot,\tk))$ & terminal cost \\ \midrule
vanilla   & 0.5998  & 213.6033    & 0.0312   \\
non-local & 1.2429  & 0.6678      & 0.1710   \\
\bottomrule
\end{tabular}
\label{tab: kitten cost compare}
\end{table}

\subsection{Computation time and accuracy}
\label{subsec: comptime}

At the end of this numerical section, we report the computation time and accuracy of above experiments in Table \ref{tab: time and acc}.
As the computational complexity of our algorithm depends on the number of vertices $\numvtc$ and number of triangles $\numtrg$ on the mesh, we also include $\numvtc,\numtrg$ in the table. 

As we mentioned in section \ref{sec: alg}, the proximal descent step requires solving a linear system \eqref{eq: mfmfg phi equ}. 
And since the linear solver is invariant to iteration numbers, we precompute it to reduce the total cost in the main iteration. The times reported in Table \ref{tab: time and acc} include both the precomputation and the main iteration.

To show that our numerical result is close to the local minimizer of the fully discretized problem \eqref{eq: mfmfg opt disct}, we report the KKT residue in Table \ref{tab: time and acc}. 
To compute the KKT residue for a given output $(\Rho,M)$, we first solve for $\Psi=\{\Psi(\cdot,\tkm)\}_{k=1,\cdots,n}\in(\prob(\disctmfd))^n$ such that
\begin{equation}
    \Dt\Dt^*\Psi - \divg\grad \Psi = 
    \Dt(\areavtc^{-1}\partial_{\Rho}\disctCost(\Rho,M)) + \divg (\areatrg^{-1}\partial_M\disctCost(\Rho,M)).
\end{equation}
Then let
\begin{equation}\left\{
\begin{aligned}
    &E_{\Rho}(\vtci,\tk):=\min\left\{\frac{1}{\areavtci}\nabla_{\Rho(\vtci,\tk)}\disctCost(\Rho,M)-(\Dt^*\Psi)(\vtci,\tk),\Rho(\vtci,\tk)\right\},\forall\vtci\in\vtc,k=1,\cdots,n\\
    &E_{M}(\trgj,\tk):=\frac{1}{\areatrgj}\nabla_{M(\trgj,\tkm)}\disctCost(\Rho,M)+(\grad\Psi)(\trgj,\tkm)=0,\forall\trgj\in\trg,k=1,\cdots,n\\
    &E_c(\vtci,\tkm):=(\Dt\Rho)(\vtci,\tkm) + (\divg M)(\vtci,\tkm) ,\forall\vtci\in\vtc,k=1,\cdots,n.
\end{aligned}\right.
\end{equation}
The KKT residue is defined as $\min\{\|E_{\Rho}\|_{\vtc,t},\|E_M\|_{\trg,t},\|E_c\|_{\vtc,t}\}$.
And $\Rho,M$ is the local minimizer of \eqref{eq: mfmfg opt disct}, if and only if the KKT residue is 0.

\begin{table}[h]
\centering
\caption{Computation time and accuracy}
\begin{tabular}{@{}cc|llllll@{}}
\toprule
\multicolumn{2}{c|}{\makecell[c]{Triangle mesh \\and interaction cost}}         & $\numvtc$                & $\numtrg$               & time(s)     & \makecell[c]{number of\\ iteration}                & \makecell[c]{time(s)\\ per iteration} &\makecell[c]{KKT\\ residue} \\ \midrule
\multicolumn{2}{c|}{U.S. map}                    & 2900                  & 5431                  & 199.4813 & 3000                 & 0.0665    & 3.15e-02         \\ \hline 
\multicolumn{2}{c|}{``8''-shape}                 & 766                   & 1536                  &  97.5769 & 5000                 & 0.0195    & 1.70e-01 \\ \hline 
\multirow{2}{*}{\makecell[c]{Irregular\\ Euclidean}} & vanilla    & \multirow{2}{*}{2473} & \multirow{2}{*}{4627} & 250.6317 & \multirow{2}{*}{5000} & 0.0501    & 1.49e-01 \\
                                     & disperse   &                       &                       & 274.1786 &                       & 0.0548    & 1.42e-01 \\ \hline 
\multirow{2}{*}{Homer}               & vanilla    & \multirow{2}{*}{2353} & \multirow{2}{*}{4702} & 163.2711 & \multirow{2}{*}{3000} & 0.0544    & 2.62e-03 \\
                                     & congested  &                       &                       & 168.6198 &                       & 0.0562    & 2.75e-03 \\ \hline 
\multirow{2}{*}{Unit shpere}         & vanilla    & \multirow{2}{*}{2562} & \multirow{2}{*}{5120} & 132.0967 & \multirow{2}{*}{2000} & 0.0660    & 6.27e-03 \\
                                     & non-local  &                       &                       & 147.9476 &                       & 0.0740    & 6.32e-03 \\ \hline 
\multirow{2}{*}{Kitten}              & vanilla    & \multirow{2}{*}{2884} & \multirow{2}{*}{5768} & 213.8950 & \multirow{2}{*}{3000} & 0.0713    & 2.99e-02 \\
                                     & non-local  &                       &                       & 260.6147 &                       & 0.0869    & 1.10e-01 \\ \bottomrule
\end{tabular}
\label{tab: time and acc}
\end{table}

\section{Conclusion}
\label{sec: conclu}

In this work, we generalize mean-field games from Euclidean space to manifolds, design an optimization based algorithm to solve variational mean field games and conduct numerical experiments on various manifolds with triangular mesh representation. 
We first propose both the PDE formulation and the variational formulation of the Nash Equilibrium of a mean-field game.
We also establish their equivalence on manifolds. 
To solve the potential mean-field games on manifolds, we use triangular meshes, piece-wise linear functions and piece-wise constant vector fields for discretization. 
Then we apply proximal gradient method to solve the corresponding discrete optimization problems.
We conduct comprehensive numerical experiments to demonstrate flexibility of the model on handling different MFG problems on various manifolds.

\bibliographystyle{plain}
\bibliography{main}

\begin{thebibliography}{10}

\bibitem{achdou2014partial}
Yves Achdou, Francisco~J Buera, Jean-Michel Lasry, Pierre-Louis Lions, and
  Benjamin Moll.
\newblock Partial differential equation models in macroeconomics.
\newblock {\em Philosophical Transactions of the Royal Society A: Mathematical,
  Physical and Engineering Sciences}, 372(2028):20130397, 2014.

\bibitem{achdou2013mean}
Yves Achdou, Fabio Camilli, and Italo Capuzzo-Dolcetta.
\newblock Mean field games: convergence of a finite difference method.
\newblock {\em SIAM Journal on Numerical Analysis}, 51(5):2585--2612, 2013.

\bibitem{achdou2010mean}
Yves Achdou and Italo Capuzzo-Dolcetta.
\newblock Mean field games: numerical methods.
\newblock {\em SIAM Journal on Numerical Analysis}, 48(3):1136--1162, 2010.

\bibitem{achdou2020mean}
Yves Achdou and Mathieu Lauri{\`e}re.
\newblock Mean field games and applications: Numerical aspects.
\newblock {\em arXiv preprint arXiv:2003.04444}, 2020.

\bibitem{almulla2017two}
Noha Almulla, Rita Ferreira, and Diogo Gomes.
\newblock Two numerical approaches to stationary mean-field games.
\newblock {\em Dynamic Games and Applications}, 7(4):657--682, 2017.

\bibitem{bauschke2011fixed}
Heinz~H Bauschke, Regina~S Burachik, Patrick~L Combettes, Veit Elser, D~Russell
  Luke, and Henry Wolkowicz.
\newblock {\em Fixed-point algorithms for inverse problems in science and
  engineering}, volume~49.
\newblock Springer Science \& Business Media, 2011.

\bibitem{beck2009fast}
Amir Beck and Marc Teboulle.
\newblock A fast iterative shrinkage-thresholding algorithm for linear inverse
  problems.
\newblock {\em SIAM journal on imaging sciences}, 2(1):183--202, 2009.

\bibitem{benamou2015augmented}
Jean-David Benamou and Guillaume Carlier.
\newblock Augmented lagrangian methods for transport optimization, mean field
  games and degenerate elliptic equations.
\newblock {\em Journal of Optimization Theory and Applications}, 167(1):1--26,
  2015.

\bibitem{benamou2017variational}
Jean-David Benamou, Guillaume Carlier, and Filippo Santambrogio.
\newblock Variational mean field games.
\newblock In {\em Active Particles, Volume 1}, pages 141--171. Springer, 2017.

\bibitem{briani2018stable}
Ariela Briani and Pierre Cardaliaguet.
\newblock Stable solutions in potential mean field game systems.
\newblock {\em Nonlinear Differential Equations and Applications NoDEA},
  25(1):1--26, 2018.

\bibitem{briceno2019implementation}
Luis Brice{\~n}o-Arias, Dante Kalise, Ziad Kobeissi, Mathieu Lauri{\`e}re,
  A~Mateos Gonz{\'a}lez, and Francisco~J Silva.
\newblock On the implementation of a primal-dual algorithm for second order
  time-dependent mean field games with local couplings.
\newblock {\em ESAIM: Proceedings and Surveys}, 65:330--348, 2019.

\bibitem{caines2018graphon}
Peter~E Caines and Minyi Huang.
\newblock Graphon mean field games and the gmfg equations.
\newblock In {\em 2018 IEEE Conference on Decision and Control (CDC)}, pages
  4129--4134. IEEE, 2018.

\bibitem{cardaliaguet2015second}
Pierre Cardaliaguet, P~Jameson Graber, Alessio Porretta, and Daniela Tonon.
\newblock Second order mean field games with degenerate diffusion and local
  coupling.
\newblock {\em Nonlinear Differential Equations and Applications NoDEA},
  22(5):1287--1317, 2015.

\bibitem{carmona2004nash}
Guilherme Carmona.
\newblock Nash equilibria of games with a continuum of players.
\newblock 2004.

\bibitem{carmona2019convergence}
Ren{\'e} Carmona and Mathieu Lauri{\`e}re.
\newblock Convergence analysis of machine learning algorithms for the numerical
  solution of mean field control and games: Ii--the finite horizon case.
\newblock {\em arXiv preprint arXiv:1908.01613}, 2019.

\bibitem{carmona2019model}
Ren{\'e} Carmona, Mathieu Lauri{\`e}re, and Zongjun Tan.
\newblock Model-free mean-field reinforcement learning: mean-field mdp and
  mean-field q-learning.
\newblock {\em arXiv preprint arXiv:1910.12802}, 2019.

\bibitem{cayton2005algorithms}
Lawrence Cayton.
\newblock Algorithms for manifold learning.
\newblock {\em Univ. of California at San Diego Tech. Rep}, 12(1-17):1, 2005.

\bibitem{de2019mean}
Antonio De~Paola, Vincenzo Trovato, David Angeli, and Goran Strbac.
\newblock A mean field game approach for distributed control of thermostatic
  loads acting in simultaneous energy-frequency response markets.
\newblock {\em IEEE Transactions on Smart Grid}, 10(6):5987--5999, 2019.

\bibitem{elie2020convergence}
Romuald Elie, Julien Perolat, Mathieu Lauri{\`e}re, Matthieu Geist, and Olivier
  Pietquin.
\newblock On the convergence of model free learning in mean field games.
\newblock In {\em Proceedings of the AAAI Conference on Artificial
  Intelligence}, volume~34, pages 7143--7150, 2020.

\bibitem{fefferman2016testing}
Charles Fefferman, Sanjoy Mitter, and Hariharan Narayanan.
\newblock Testing the manifold hypothesis.
\newblock {\em Journal of the American Mathematical Society}, 29(4):983--1049,
  2016.

\bibitem{GLM}
Wilfrid Gangbo, Wuchen Li, and Chenchen Mou.
\newblock Geodesics of minimal length in the set of probability measures on
  graphs.
\newblock {\em ESAIM: COCV}, 2019.

\bibitem{gao2021modeling}
Hao Gao, Wuchen Li, Miao Pan, Zhu Han, and H~Vincent Poor.
\newblock Modeling covid-19 with mean field evolutionary dynamics: Social
  distancing and seasonality.
\newblock {\em Journal of Communications and Networks}, 23(5):314--325, 2021.

\bibitem{gomes2018mean}
Diogo Gomes and Jo{\~a}o Sa{\'u}de.
\newblock A mean-field game approach to price formation in electricity markets.
\newblock {\em arXiv preprint arXiv:1807.07088}, 2018.

\bibitem{gomes2013continuous}
Diogo~A Gomes, Joana Mohr, and Rafael~Rigao Souza.
\newblock Continuous time finite state mean field games.
\newblock {\em Applied Mathematics \& Optimization}, 68(1):99--143, 2013.

\bibitem{gomes2021numerical}
Diogo~A Gomes and Jo{\~a}o Sa{\'u}de.
\newblock Numerical methods for finite-state mean-field games satisfying a
  monotonicity condition.
\newblock {\em Applied Mathematics \& Optimization}, 83(1):51--82, 2021.

\bibitem{gueant2015existence}
Olivier Gu{\'e}ant.
\newblock Existence and uniqueness result for mean field games with congestion
  effect on graphs.
\newblock {\em Applied Mathematics \& Optimization}, 72(2):291--303, 2015.

\bibitem{huang2007large}
Minyi Huang, Peter~E Caines, and Roland~P Malham{\'e}.
\newblock Large-population cost-coupled lqg problems with nonuniform agents:
  individual-mass behavior and decentralized $\varepsilon$-nash equilibria.
\newblock {\em IEEE transactions on automatic control}, 52(9):1560--1571, 2007.

\bibitem{huang2006large}
Minyi Huang, Roland~P Malham{\'e}, Peter~E Caines, et~al.
\newblock Large population stochastic dynamic games: closed-loop mckean-vlasov
  systems and the nash certainty equivalence principle.
\newblock {\em Communications in Information \& Systems}, 6(3):221--252, 2006.

\bibitem{kuhn1951nonlinear}
H.~W. Kuhn and A.~W. Tucker.
\newblock Nonlinear programming.
\newblock In {\em Proceedings of the Second Berkeley Symposium on Mathematical
  Statistics and Probability}, pages 481--492. University of California Press,
  1951.

\bibitem{lai2011framework}
Rongjie Lai and Tony~F Chan.
\newblock A framework for intrinsic image processing on surfaces.
\newblock {\em Computer vision and image understanding}, 115(12):1647--1661,
  2011.

\bibitem{lasry2007mean}
Jean-Michel Lasry and Pierre-Louis Lions.
\newblock Mean field games.
\newblock {\em Japanese journal of mathematics}, 2(1):229--260, 2007.

\bibitem{lauriere2021numerical}
Mathieu Lauri{\`e}re.
\newblock Numerical methods for mean field games and mean field type control.
\newblock {\em arXiv preprint arXiv:2106.06231}, 2021.

\bibitem{lee2013smooth}
John~M Lee.
\newblock Smooth manifolds.
\newblock In {\em Introduction to Smooth Manifolds}, pages 1--31. Springer,
  2013.

\bibitem{lee2017global}
Taeyoung Lee, Melvin Leok, and N~Harris McClamroch.
\newblock Global formulations of lagrangian and hamiltonian dynamics on
  manifolds.
\newblock {\em Springer}, 13:31, 2017.

\bibitem{LLO1}
Wonjun Lee, Wuchen Li, and Stanley Osher.
\newblock Mean field control problems for vaccine distribution.
\newblock {\em arXiv:2104.11887}, 2021.

\bibitem{LLTLO}
Wonjun Lee, Siting Liu, Hamidou Tembine, Wuchen Li, and Stanley Osher.
\newblock Controlling propagation of epidemics via mean-field control.
\newblock {\em SIAM Journal on Applied Mathematics}, 81(1):190--207, 2021.

\bibitem{LiLeeSO}
Wuchen Li, Wonjun Lee, and Stanley Osher.
\newblock Computational mean-field information dynamics associated with
  reaction diffusion equations.
\newblock {\em arXiv:2107.11501}, 2021.

\bibitem{LiLiuSO}
Wuchen Li, Siting Liu, and Stanley Osher.
\newblock Controlling conservation laws i: entropy-entropy flux.
\newblock {\em arXiv:2111.05473}, 2021.

\bibitem{lin2020apac}
Alex~Tong Lin, Samy~Wu Fung, Wuchen Li, Levon Nurbekyan, and Stanley~J Osher.
\newblock Apac-net: Alternating the population and agent control via two neural
  networks to solve high-dimensional stochastic mean field games.
\newblock {\em PNAS}, 2021.

\bibitem{liu2021computational}
Siting Liu, Matthew Jacobs, Wuchen Li, Levon Nurbekyan, and Stanley~J Osher.
\newblock Computational methods for first-order nonlocal mean field games with
  applications.
\newblock {\em SIAM Journal on Numerical Analysis}, 59(5):2639--2668, 2021.

\bibitem{MAAS20112250}
Jan Maas.
\newblock Gradient flows of the entropy for finite markov chains.
\newblock {\em Journal of Functional Analysis}, 261(8):2250--2292, 2011.

\bibitem{mangasarian1975pseudo}
Olvi~L Mangasarian.
\newblock Pseudo-convex functions.
\newblock In {\em Stochastic optimization models in finance}, pages 23--32.
  Elsevier, 1975.

\bibitem{meyer2003discrete}
Mark Meyer, Mathieu Desbrun, Peter Schr{\"o}der, and Alan~H Barr.
\newblock Discrete differential-geometry operators for triangulated
  2-manifolds.
\newblock In {\em Visualization and mathematics III}, pages 35--57. Springer,
  2003.

\bibitem{nash1951non}
John Nash.
\newblock Non-cooperative games.
\newblock {\em Annals of mathematics}, pages 286--295, 1951.

\bibitem{nurbekyan2019fourier}
Levon Nurbekyan and Sa{\'u}de Jo{\~a}o.
\newblock Fourier approximation methods for first-order nonlocal mean-field
  games.
\newblock {\em Portugaliae Mathematica}, 75(3):367--396, 2019.

\bibitem{rockafellar1970convex}
R~Tyrrell Rockafellar.
\newblock {\em Convex analysis}, volume~36.
\newblock Princeton university press, 1970.

\bibitem{ruthotto2020machine}
Lars Ruthotto, Stanley~J Osher, Wuchen Li, Levon Nurbekyan, and Samy~Wu Fung.
\newblock A machine learning framework for solving high-dimensional mean field
  game and mean field control problems.
\newblock {\em Proceedings of the National Academy of Sciences},
  117(17):9183--9193, 2020.

\bibitem{solomon2016entropic}
Justin Solomon, Gabriel Peyr{\'e}, Vladimir~G Kim, and Suvrit Sra.
\newblock Entropic metric alignment for correspondence problems.
\newblock {\em ACM Transactions on Graphics (TOG)}, 35(4):1--13, 2016.

\bibitem{weinan2019mean}
E~Weinan, Jiequn Han, and Qianxiao Li.
\newblock A mean-field optimal control formulation of deep learning.
\newblock {\em Research in the Mathematical Sciences}, 6(1):10, 2019.

\bibitem{yang2017mean}
Chungang Yang, Jiandong Li, Min Sheng, Alagan Anpalagan, and Jia Xiao.
\newblock Mean field game-theoretic framework for interference and energy-aware
  control in 5g ultra-dense networks.
\newblock {\em IEEE Wireless Communications}, 25(1):114--121, 2017.

\bibitem{yang2018mean}
Yaodong Yang, Rui Luo, Minne Li, Ming Zhou, Weinan Zhang, and Jun Wang.
\newblock Mean field multi-agent reinforcement learning.
\newblock In {\em International Conference on Machine Learning}, pages
  5571--5580. PMLR, 2018.

\bibitem{yu2021fast}
Jiajia Yu, Rongjie Lai, Wuchen Li, and Stanley Osher.
\newblock A fast proximal gradient method and convergence analysis for dynamic
  mean field planning.
\newblock {\em arXiv preprint arXiv:2102.13260}, 2021.

\end{thebibliography}

\end{document}